\documentclass{amsart}
\usepackage{amsmath,amssymb,amscd,amsthm,amsxtra,enumerate,amsfonts}
\usepackage{latexsym}

\newtheorem{theorem}{Theorem}
\newtheorem{lemma}{Lemma}

\newtheorem{claim}{Claim}
\newtheorem{case}{Case}
\newtheorem*{subcase4a}{Subcase 4(a)}
\newtheorem*{subcase4b}{Subcase 4(b)}

\newcommand{\q}{\quad}
\newcommand{\qq}{\quad\quad}

\newcommand{\De}{\Delta}

\newcommand{\om}{\omega}
\newcommand{\Om}{\omega}
\newcommand{\oms}{\omega_p}
\newcommand{\omu}{\omega_u}
\newcommand{\omsr}{\omega_{p(r)}}
\newcommand{\omsi}{\omega_{p(i)}}
\newcommand{\omsone}{\omega_{p(1)}}
\newcommand{\is}{I_p}
\newcommand{\iuu}{I_u}
\newcommand{\itt}{I_{\tr}}
\newcommand{\itr}{I_{\tr'}}
\newcommand{\itj}{I_{\tr_j}}
\newcommand{\x}{\mathbf x}

\newcommand{\cm}{\mathcal M}

\newcommand{\ce}{\mathcal E}

\newcommand{\dd}{\mathbf D}
\newcommand{\pp}{\mathbf P}
\newcommand{\tr}{\mathbf T}
\newcommand{\bu}{\mathbf U}
\newcommand{\bv}{\mathbf V}
\newcommand{\N}{\mathbf N}

\newcommand{\lab}{\label}

\newcommand{\intl}{\int\limits}

\newcommand{\intrn}{\int_{\rn}}

\newcommand{\dint}{\displaystyle\int}

\newcommand{\f}{\frac}
\newcommand{\df}{\displaystyle\frac}

\newcommand{\wh}{\widehat}

\newcommand{\nf}{\infty}

\newcommand{\rn}{\mathbb R^n}

\newcommand{\zzz}{\mathbf Z}
\newcommand{\lng}{\langle}
\newcommand{\rng}{\rangle}
\newcommand{\summ}{\displaystyle\sum}
\newcommand{\supp}{\displaystyle\sup}

\begin{document}

\title
[A weak $L^2$ estimate for a maximal dyadic sum operator on $\mathbb R^n$]
{A weak $L^2$ estimate for a maximal dyadic sum operator on $\mathbb R^n$}

\author{Malabika Pramanik}

\address{Malabika Pramanik \\ Department of Mathematics \\ University of Wisconsin, Madison \\ Madison, WI 53706}

\email{malabika@math.wisc.edu}

\author{Erin Terwilleger}

\address{
Erin Terwilleger\\
School of Mathematics\\
Georgia Institute of Technology\\
Atlanta, GA 30332 USA}

\email{erin@math.gatech.edu}

\date{\today}



\begin{abstract}
Lacey and Thiele have recently obtained a new proof of Carleson's theorem on almost everywhere convergence
of Fourier series. This paper is a generalization of their techniques (known broadly as time-frequency
analysis) to higher dimensions. In particular, a weak-type (2,2) estimate is derived for a maximal dyadic
sum operator on $\mathbb R^{n}$, $n > 1$. As an application one obtains a new proof of Sj\"olin's theorem
on weak $L^{2}$ estimates for the maximal conjugated Calder\'on-Zygmund operator on $\mathbb R^{n}$.  
\end{abstract}
\maketitle

\section{Introduction}

In 1966, Carleson \cite{carleson} proved his celebrated theorem on almost everywhere convergence of Fourier
series of square integrable functions on $\mathbb R$.  This was followed by a new proof 
given by C. Fefferman \cite{cfefferman} in 1973.  The techniques used by C. Fefferman
 have become known as time-frequency analysis and
have found wide application in harmonic analysis in recent years.  In particular,
Lacey and Thiele \cite{lacey-thiele1}, \cite{lacey-thiele2} have refined and extended
these ideas in their pioneering work on the
 bilinear Hilbert transform on $\mathbb R$.  In 2000, they obtained a new proof of Carleson's theorem
 \cite{lacey-thiele3} in which these
 techniques play a crucial role.  These powerful techniques stem from interaction of extremely deep
 ideas which include delicate orthogonality estimates, combinatorics, and quasi-orthogonal
 decompositions well-localized in both time and space.  It is the goal of this paper to extend
 the techniques of time-frequency analysis of \cite{lacey-thiele3} to higher dimensions.

The main result of this paper is a weak-type $L^{2}$ estimate for a maximal dyadic
sum operator in $\mathbb R^{n}$, $n > 1$. 
In dimension one, this operator may be thought of as a
linearized and discretized version of the Carleson operator $\mathcal C$,
\[ \mathcal Cf(x) := \sup_{N} \left| \int_{- \infty}^{N} \widehat{f}(\xi) e^{2 \pi i x \xi} \, d\xi
\right|.\] The main point of Lacey and Thiele's proof \cite{lacey-thiele3} is to show that the 
discretized
operator satisfies a weak-type (2,2) estimate. This in turn implies a similar
estimate for $\mathcal C$, which is a key ingredient in proving that the Fourier
series of a square-integrable function on the circle converges almost everywhere.

We introduce a higher dimensional analogue of the linearized and
discretized Carleson operator and adapt the methodology of Lacey and Thiele to
prove that this operator maps $L^{2}(\mathbb R^{n})$ to $L^{2, \infty}(\mathbb R^{n})$.
One of the distinguishing aspects of our proof is the introduction of an ordering
of points in $\rn$ which allows us to organize the higher dimensional rectangles
 and thus control the large sums that appear in the operator.
Unlike the situation in dimension one, the mapping property above does not lead to an almost
everywhere convergence result in higher dimensions. However, it gives as a corollary a
result of Sj\"olin \cite{sjolin} on the weak $L^{2}$ boundedness of the maximal conjugated
Calder\'on-Zygmund operator on $\mathbb R^{n}$.

The proof is divided into seven sections. The first section explains the notation
and terminology and gives the statement of the main theorem. The second section
lists the main ingredients of the proof and the argument that binds them together.
The subsequent four sections are devoted to proving the different lemmas needed in 
the main argument. The final section provides a new proof
of Sj\"olin's theorem as an application of our main result.

\section{Main Theorem}
\setcounter{equation}0 

Time-frequency analysis provides the  crucial 
set of ideas in the recent progress made in the understanding of 
Carleson's theorem. In this type of analysis one
heavily uses the structure of dyadic intervals.  A dyadic interval has the 
form $[m2^{k},(m+1)2^{k})$ where $k$ and $m$ are integers and $k$ is called 
the scale.  A dyadic cube $I \subset \rn$ 
 is of the form
 $$
 \prod_{j=1}^n I^j = \prod_{j=1}^n \, \, [m_j 2^k,(m_j +1) 2^k),
 $$
 where  $k$ and $m_j$ are integers for all $j= 1,2, \cdots, n$. 
 We easily see that the $n$-dimensional volume
is given by $|I|=2^{nk}$.
Let  $c(I)=(c(I^1), \cdots , c(I^n))$ denote the center of $I$, and for   
$a>0$, $aI$ will denote the cube with the same center as $I$ and whose volume
is $a^n|I|$.

Consider the time-frequency plane in $2n$ dimensional space with points $(\x,{\pmb {\xi}})$, where 
$\x$ 
denotes the time coordinate in $\rn$ and ${\pmb {\xi}}$ denotes the frequency coordinate in $\rn$.  A 
``rectangle'' in the time-frequency plane is the cross product of a dyadic
cube from the time plane and a dyadic cube from the frequency plane.  To be more precise, for 
a rectangle, $p$, the projection onto the time plane will be denoted by $I_p$, and its projection
onto the frequency plane will be denoted by $\om_p$.  We will denote by $\dd$ the set of rectangles
$p=I_p \times \om_p$ such that $|I_p||\oms|=1$.  An element of $\dd$ will be called a {\em {tile}}.

As mentioned earlier,  it is important for the higher dimensional version of our
time-frequency analysis to introduce an ordering in $\rn$ that will play a role analogous
to the linear ordering on $\mathbb R$.  This is especially relevant in a certain selection scheme used
in Section 5 in analogy with the work of Lacey and Thiele.  Although the choice of ordering is not
unique (we will mention an alternative in Section 5),
 we find it convenient to work with the lexicographical order defined as follows.
 Given $\mathbf a =(a_1, a_2, \cdots, a_n), \mathbf b =(b_1, b_2, \cdots, b_n)
\in \rn$,

$$
\mathbf a < \mathbf b \iff
\begin {cases}
    a_1 < b_1 &\\
    a_1 =  b_1 \, , \, a_2 < b_2 &\\
    \q \vdots &\\
  a_1 =  b_1\, , \, a_2 = b_2\, , \,  \ldots \, , \, a_{n-1} = b_{n-1} \, , \, a_n < b_n \, , &
\end{cases}
$$
where the right hand side above is to be read with Boolean ``or'' standard.

For a tile $p$ with $\om_p=\om_p^1 \times \om_p^2 \times \ldots
\times \om_p^n $, we can divide each dyadic interval $\om_p^j$ into two parts.  In
other words, for $j=1,2, \ldots, n$, we get 
$$
\om_p^j= (\om_p^j \cap (-\nf,c(\om_p^j)) \cup 
(\om_p^j \cap [c(\om_p^j),\nf)).
$$
 Then $\om_p$ can be decomposed into 
$2^n$ subcubes formed from all combinations
of cross products of these half intervals.  We number these subcubes using the
lexicographical order on the centers and denote the subcubes by $\omsi$
for $i=1,2, \ldots, 2^n$.  A tile $p$ is then the union of $2^n$
{\em {semi-tiles}} given by  
$p(i)=I_p \times \omsi$ for $i=1,2, \ldots, 2^n$. 

Let us define translation, modulation, and dilation operators by
\begin{align*}
& T_{\bf y} f(\x) := f(\x-\bf y) \\
& M_{{\pmb {\eta}}} f(\x) := f(\x) e^{2 \pi i {\pmb {\eta}} \cdot \x} \\
& D_{\lambda}^q f(\x) := \lambda^{-n/q} f(\lambda^{-1}\x), \lambda > 0.
\end{align*}
Note that if we set $q=2$, these operators are isometries on $L^2(\rn)$.
We fix a Schwartz function $\phi$ such that
$\widehat {\phi}$ is real, nonnegative, supported in the cube
$[-1/10,1/10]^n$ and equal to 1 on the cube 
$[-9/100,9/100]^n$. For a tile $p \in \dd$ and $\x \in \rn$ we define 
\begin{equation}\lab{31}
\phi_p (\x)=M_{c(\omsone)} T_{c(I_p)} D_{|I_p|^ { 1/n}}^2 \phi (\x).
\end{equation}
Using the following definition of  the Fourier transform  
$$
\wh{f}({\pmb {\xi}}) = \int_{\rn} f(x) e^{-2 \pi i \x \cdot {\pmb {\xi}}}\, d\x, 
$$
one easily can see that
\begin{equation}\lab{32}
{\widehat {\phi_p}}({\pmb {\xi}})= T_{c(\omsone)} M_{-c(I_p)} D_{|\oms|^{1/n}}^2 
\widehat {\phi}({\pmb {\xi}}).
\end{equation}
Equation (\ref{31}) tells us that for each $p$ the function $\phi_p$ is well 
localized in time with most of its mass in $I_p$ while
equation (\ref{32}) tells us that $\widehat {\phi_p}$ is supported 
in $\frac {1} {5} \omsone$.  Note also that the $\phi_p$ have the same $L^2(\rn)$ norm.

Let $m$ be a multiplier  in $C^{\infty}(\rn \setminus \{\mathbf 0\})$ which is homogeneous of degree 0,
and define
$$
  \left (\psi_p^{{\pmb {\zeta}}}\right)\widehat \, \, ({\pmb{\xi}}) = m({\pmb {\xi}} - {\pmb {\zeta}})
\widehat {\phi_p}({\pmb {\xi}}),
$$
where ${\pmb {\zeta}}$ is contained in $\omsr$ for some fixed $r \in \{2,3, \ldots, 2^n\}$.
Note that we have the following fact for all ${\pmb {\zeta}} \in \omsr$:
$$
|\psi_p^{{\pmb {\zeta}}}(\x)| \le C_{\nu} |I_p|^{-1/2} \Bigl( 1 + \f {|\x -c(I_p)|} { |I_p|^{1/n}} \Bigr)^{-\nu},
$$
where $\nu$ is a large integer whose value may vary at different places in the proof.
To see this fact we write
\begin{align*}
\psi_p^{{\pmb {\zeta}}}(\x) = & \int_{\rn} e^{2 \pi i {\pmb {\xi}} \cdot \x} \widehat 
{\phi_p}({\pmb {\xi}})
m({\pmb {\xi}} -{\pmb {\zeta}}) \, d{\pmb {\xi}} \\
= &  \int_{\rn} e^{2 \pi i {\pmb {\xi}} \cdot \x} T_{c(\omsone)} M_{-c(I_p)} D_{|\oms|^{1/n}}^2 
\widehat {\phi}({\pmb {\xi}}) m({\pmb {\xi}} -{\pmb {\zeta}}) \, d{\pmb {\xi}} \\
= & e^{2 \pi i c(I_p) \cdot c(\om_{p(1)})}  \\
& \times \int_{\rn} e^{2 \pi i {\pmb {\xi}} \cdot (\x- c(I_p))}
|I_p|^{1/2} \widehat \phi \bigl( |I_p|^{1/n}({\pmb {\xi}}-c(\om_{p(1)})) \bigr) 
m({\pmb {\xi}} -{\pmb {\zeta}}) \, d{\pmb {\xi}}. 
\end{align*}
Making the change of variable
$$
 {\pmb {\xi'}} = |I_p|^{1/n} ({\pmb {\xi}} - c(\om_{p(1)}))
$$
in the above integral, we obtain
\begin{align*}
\psi_p^{{\pmb {\zeta}}}(\x) = & 
e^{2 \pi i c(I_p) \cdot c(\om_{p(1)})} \int_{\rn} e^{2 \pi i \bigl({\pmb {\xi'}}|I_p|^{-1/n} + 
c(\om_{p(1)}) \bigr)
\cdot (\x- c(I_p))}
|I_p|^{-1/2} \widehat \phi ({\pmb {\xi'}}) \\ 
& \hskip 5cm \times m \bigl({\pmb {\xi'}} |I_p|^{-1/n} + c(\om_{p(1)})-{\pmb {\zeta}} \bigr) \, d{\pmb {\xi'}} \\
= & e^{2 \pi i \x \cdot c(\om_{p(1)})} \int_{\rn} e^{2 \pi i {\pmb {\xi'}} \cdot 
(\x- c(I_p))|I_p|^{-1/n}} 
|I_p|^{-1/2} \widehat \phi ({\pmb {\xi'}})  \\
& \hskip 5cm \times m \bigl({\pmb {\xi'}} + |I_p|^{1/n}(c(\om_{p(1)})-{\pmb {\zeta}}) \bigr) \, d{\pmb {\xi'}},
\end{align*}
where the second equality is obtained using the fact that $m$ is homogeneous of degree 0.
Since ${\pmb {\zeta}} \in \om_{p(r)}$ and ${\pmb {\xi}} \in \f {1} {2} \om_{p(1)}$,
we have
$$
|{\pmb {\zeta}}-{\pmb {\xi}}| \gtrsim |\om_p|^{1/n},
$$
from which it follows that
$$
{\pmb {\xi'}} + |I_p|^{1/n}(c(\om_{p(1)})-{\pmb {\zeta}})) \gtrsim 1.
$$
Thus all the derivatives of 
$m({\pmb {\xi'}} + |I_p|^{1/n}(c(\om_{p(1)})-{\pmb {\zeta}}))$ are bounded.  A standard
integration by parts argument finishes the proof.

Using the following definition for the inner product
$$
\lng f,g \rng = \intrn f(\x) \overline {g(\x)} \, d\x \, ,
$$
given ${\pmb {\zeta}} \in \rn$ and $f \in L^{2}(\rn)$,  
we define an operator
$$
B_{{\pmb {\zeta}}}^r f(\cdot)=\sum _{p \in \dd}   \langle
f,\phi_p \rangle \psi_p^{{\pmb {\zeta}}}(\cdot) 1_{\omsr}({\pmb {\zeta}}).
$$
\begin{theorem}\label{thm}
  There exists a constant C, depending only on dimension,
 so that
for all $f \in L^2(\rn)$ and $r \in \{2,3, \ldots, 2^n\}$
\begin{equation}\lab{33}
\| \displaystyle\sup _{{\pmb {\zeta}}\in \rn }|B_{{\pmb {\zeta}}}^r f|\|_{L^{2,\nf}(\rn)} 
\le C \|f\|_{L^2(\rn)}.
\end{equation}
\end{theorem}

To prove the theorem, we will work with a linearized version of the operator.  Consider a 
measurable function $\x \rightarrow \N(\x)=(N_1(\x),N_2(\x), \ldots, N_n(\x))$ 
from $\rn$ to $\rn$ and define a linear operator
$$
B_\N^{r}(\x):=  B_{\N(\x)}^{r}(\x)= \sum_{p \in \dd}  \lng f,\phi_p \rng \psi_p^{\N(\x)}(\x) 
(1_{\omsr} \circ \N)(\x).
$$
To prove (\ref{33}) it will suffice to show that there exists a constant $C>0$ such that for 
all $f \in L^2(\rn)$ 
\begin{equation}\lab{34}
\displaystyle\sup _{\N:\rn \rightarrow \rn} \|B_{\N}^r f \|_{L^{2, \nf}(\rn)} \le 
C \|f\|_{L^2(\rn)},
\end{equation}
where the supremum is taken over all measurable functions $\N$ on $\rn$.

By duality we will show that the adjoint operator
$$
g \rightarrow \summ_{p \in \dd} \lng (1_{\omsr} \circ \N) \psi_p^{\N},g \rng \phi_p
$$
maps $L^{2,1}(\rn)$ into $L^2(\rn)$ 
with bounds independent of the measurable function $\N$. 
Since $L^{2,1}(\rn)$ is a Lorentz space, it suffices to show that the dual operator maps 
$L^{2,1}(\rn) \cap \{1_E : E \subset \rn, E {\mbox{ measurable }}, |E| < \infty \}$ into
$L^2(\rn)$. Hence, we need to show
\begin{equation}\lab{36}
\| \summ_{p \in \dd} \lng (1_{\omsr} \circ \N) \psi_p^{\N},1_E \rng \phi_p \|_{L^2(\rn)}  \le
C|E|^{1/2}.
\end{equation}
By duality, (\ref{36}) is equivalent to
\begin{equation}\lab{37}
| \summ_{p \in \dd} \lng (1_{\omsr} \circ \N) \psi_p^{\N},1_E \rng \lng \phi_p, f \rng |
  \le C|E|^{1/2},
\end{equation}
for all Schwartz functions $f$ with $L^2$ norm one.  We will further restrict the sum
to an arbitrary finite subset $\mathbf P$ of $\dd$.

Now for all integers $j$ we have the identity
\begin{multline*}
 \summ_{p \in \pp}| \lng (1_{\omsr} \circ \N) \psi_p^\N,1_E \rng \lng \phi_p, f \rng |
 = \\ 2^{- \f {jn} {2}}\sum_{u \in P(j)}| \lng (1_{\om_{u(r)}} \circ \N_j) \psi_u^{\N_j},
1_{2^j \otimes E} \rng \lng \phi_u, 2^{-\f {jn} {2}}f(2^{-j}(\cdot) )\rng |,
\end{multline*}
where for any set $A$ we define $2^j \otimes A = \{2^j \mathbf y =(2^jy_1,2^jy_2, \ldots, 2^jy_n):
\mathbf y \in A\}$, $\N_j(\x)= 2^{-j} \N(2^{-j}\x)$, and $\pp(j) = \{(2^j \otimes I_p) \times 
(2^{-j} \otimes \oms):p \in \pp \}$.  By picking $j$ so that $1 \le 2^{jn}|E| \le 2$, we can 
absorb $|E|$ into the constant on the right hand side of (\ref{37}).  Finally we note that the
left hand side of (\ref{37}) can be rewritten so that the estimate we need to show now becomes
\begin{equation}\lab{39}
 \summ_{p \in \pp}| \lng 1_{E \cap \N^{-1}[\omsr]} , \psi_p^\N \rng \lng \phi_p, f \rng |
  \le C,
\end{equation}
for all Schwartz functions $f$ with $L^2$ norm one, measurable functions $\N$, measurable sets 
$E$ with $|E| \le 1$, and all finite subsets $\pp$ of $\dd$.  For the rest of the paper we
fix $f$, $\N$, and $E$ in this manner.  By $\N^{-1}[\omsr]$ we mean
$\{\x: \N(\x) \in \omsr \}$.

\vspace{5mm}
\section{Main Argument}
\setcounter{equation}0

We now set up some tools that we will use throughout the rest of the paper.  Define
a partial order $<$ on the set of tiles $\dd$ by setting
\begin{center}
$p<p' \qq \iff \qq I_p \subset I_{p'}$ \q and \q $\om_{p'} \subset \oms$.
\end{center}
We have the property that if two tiles $p,p' \in \dd$ intersect, then either $p<p'$ or $p'<p$.
To see this, observe that dyadic cubes have the property that if two of them intersect, then one
is contained in the other.  
This extends from the same property for dyadic intervals in dimension one. Now, suppose two
tiles $p$ and $p'$ in $\dd$ intersect, and without loss of generality let $|I_{p}| \leq |I_{p'}|$.
Then
$p$ and $p'$ intersect in both 
the time and frequency components, i.e. $I_{p} \cap I_{p'} \ne \emptyset, \; \om_{p} \cap
\om_{p'} \ne \emptyset$. From size considerations, one obtains
that $I_p \subset I_{p'}$ and $\om_{p'} \subset \om_{p}$, hence $p < p'$. 
A consequence of this property is that for a finite set of tiles $\pp$, all maximal 
elements of $\pp$ under $<$ must be disjoint sets.

A finite set of tiles $\tr$ is called a {\em {tree}} if there
exits a tile $t \in \dd$ such that $p<t$ for all $p \in \tr$.  We call $t$ the {\em {top}} of the tree
$\tr$ and denote it by $p_{\tr}=I_{\tr} \times \om_{\tr}$.  Note that the top is unique but not 
necessarily an
element of the tree.  Another useful observation is that any finite set of tiles $\pp$ can
be written as a union of trees. Consider all maximal elements of $\pp$ under $<$.  Then
a nonmaximal element $p \in \pp$ must be less than, under ``$<$'', some maximal element $t \in \pp$
which places $p$ in the tree with top $t$.
For $i \in \{1,2, \ldots, 2^n\}$, we call a tree an {\em {$i$-tree}}, denoted by $\tr^i$, if
$
\om_{\tr(i)} \subset \omsi
$
for all $p \in \tr^i$.  Observe that any tree can be written as the disjoint union of $i$-trees. 
Also
for fixed $i_0$, and $p,p' \in \tr^{i_0}$, the subcubes $\omsi$ and $\om_{p'(i)}$ are pairwise disjoint and disjoint from
 $\om_{\tr(i)}$ for all
$i \in \{1,2, \ldots, 2^n\} \setminus \{i_0\}$.

For $p \in \dd$, define the {\em {mass}} of $\{p\}$ as
$$
\cm(\{p\})= \supp_{\substack {u \in \dd \\ p<u }} \ \intl_{E \cap \N^{-1}[\om_u]} \df {|\iuu|^{-1}} 
{\bigl(1+ \f {|\x-c(\iuu)|} {|\iuu|^{1/n}}\bigr)^{10n}} \, d\x.
$$                   
We can then define the mass of a finite set of tiles $\pp$ to be
$$
\cm(\pp)= \supp_{p \in \pp} \cm(\{p\}).
$$
Note that the mass of any set of tiles is at most one since by a change of variables
$$
\cm(\pp) \le \dint_{\rn} \df {1} {(1 + |\x|)^{10n}} \, d\x \le 1.
$$
The {\em {energy}}, depending on $r$, of a finite set of tiles $\pp$ is defined
$$
\ce(\pp)=\supp_{\tr^r \in \pp} \Bigl(|I_{\tr^r}|^{-1} \summ_{p \in \tr^r} 
|\lng f, \phi_p \rng|^2 \Bigr)^{1/2}.
$$
Recall that $r \in \{2,3, \ldots, 2^n\}$ is fixed and $f$ is a fixed Schwartz function of $L^2(\rn)$ 
norm one.
The following three lemmata will provide the main steps in proving the theorem, and their proofs will 
be shown in the next four sections of the paper.


\begin{lemma}\lab{lem1}
There exists a constant $C_1$ such that for any finite set of tiles $\pp$ there is a subset $\pp'$ 
of $\pp$ such that
\begin{equation}\lab{41}
\cm(\pp \setminus \pp') \le \f {1} {4} \cm(\pp)
\end{equation}
and $\pp'$ is the union of trees $\tr_j$ satisfying
\begin{equation}\lab{42}
\summ_j |\itj| \le \f {C_1} {\cm(\pp)}.
\end{equation}
\end{lemma}


\begin{lemma}\lab{lem2}
There exists a constant $C_2$ such that for any finite set of tiles $\pp$ there is a subset $\pp''$ 
of $\pp$ such that
\begin{equation}\lab{43}
\ce(\pp \setminus \pp'') \le \f {1} {2} \ce(\pp)
\end{equation}
and $\pp''$ is the union of trees $\tr_j$ satisfying
\begin{equation}\lab{44}
\summ_j |\itj| \le \f {C_2} {\ce(\pp)^2}.
\end{equation}
\end{lemma}


\begin{lemma}\lab{lem3}
(The Tree Inequality)  There exists a constant $C_3$ such that for all trees $\tr$  
\begin{equation}\lab{45}
\summ_{p \in \tr}| \lng 1_{E \cap \N^{-1}[\omsr]} , \psi_p^{\N} \rng \lng \phi_p, f \rng |
  \le C_3 |\itt| \ce(\tr) \cm(\tr).
\end{equation}
\end{lemma}

We will now prove (\ref{39}), and hence Theorem \ref{thm}, assuming the three lemmata.
In the argument below set
$$
C_0=C_1 + C_2.
$$
Given a finite set of tiles $\mathbf P$, 
find   a very large integer $m_0$ such that $\mathcal E(\mathbf P) \le 2^{m_0n}$
and $\mathcal M(\mathbf P)\le 2^{2m_0n}$. We construct by 
decreasing induction a
sequence of  pairwise disjoint sets  $\mathbf P_{m_0}$, 
$\mathbf P_{m_0-1}$, $\mathbf P_{m_0-2}$,  $\mathbf P_{m_0-3}$, ... such
that 
$$
\bigcup_{j=-\nf}^{m_0} \mathbf P_{j} = \mathbf P
$$
and such that the following properties are satisfied 
\begin{enumerate}
\item[(1)] $\mathcal E(\mathbf P_{ j}) \le 2^{(j+1)n} $   for all $j\le m_0$.
\item[(2)] $\mathcal M(\mathbf P_{ j}) \le 2^{(2j+2)n}$   for all $j\le m_0$.
\item[(3)] $\mathcal E\big(\mathbf P \setminus (\mathbf P_{m_0}  
\cup \dots \cup \mathbf P_{j} )\big)  \le 2^{ jn } $   for all $j\le m_0$.
\item[(4)] $\mathcal M\big(\mathbf P \setminus (\mathbf P_{m_0} 
\cup \dots \cup\mathbf P_{j} )\big)  \le 2^{ 2jn }$   for all $j\le m_0$.
\item[(5)] $\mathbf P_j$ is a union of trees $\mathbf T_{jk}$ 
such that 
$\sum_{k} |I_{{\mathbf T_{jk}}}| \le C_0 2^{-2jn}$ 
for all $j\le m_0$.
\end{enumerate}

Assume momentarily that we have constructed a sequence 
$\mathbf P_j$ as above. Then to obtain estimate (\ref{39}) we
use (1), (2), (5), the observation that the mass is always bounded by $1$,
and    Lemma  \ref{lem3} to obtain 
{\allowdisplaybreaks
\begin{align*}
& \sum_{s\in \mathbf P}  \, \big|
  \langle     1_{E\cap \N^{-1}[\om_{p(r)}]} , \psi_p^{\N} \rangle 
\langle  f,   \phi_p  \rangle   \big| \\
\le &\sum_{j} \sum_{p\in \mathbf P_{j }}  \, \big|
  \langle     1_{E\cap \N^{-1}[\om_{p(r)}]} , \psi_p^{\N} \rangle
\langle  f,   \phi_p  \rangle   \big| \\
\le &\sum_{j}\sum_{k} \sum_{p\in \mathbf T_{jk}}  \, \big|
 \langle     1_{E\cap \N^{-1}[\om_{p(r)}]} , \psi_p^{\N} \rangle 
\langle  f,   \phi_p  \rangle  \big| \\
\le &C_3\sum_{j}\sum_{k}  
|I_{\mathbf T_{jk}}|\, \mathcal E
(\mathbf T_{jk})\, \mathcal M(\mathbf T_{jk}) \\
\le &C_3\sum_{j}\sum_{k}  
|I_{\mathbf T_{jk}}|\, 2^{(j+1)n}\, \min (1, 2^{(2j+2)n}) \\
\le & C_3 \sum_{j}  C_0   2^{-2jn} 2^{(j+1)n}\, \min (1, 2^{(2j+2)n})  \\
\le & 2^{3n} C_0C_3 \sum_{j}\min(2^{jn},2^{-jn}) \le C_n\, . 
\end{align*}}
This proves estimate (\ref{39}). 

It remains to construct a sequence of disjoint sets   $\mathbf P_j$ 
satisfying (1)-(5).  We start our induction at $j=m_0$ by 
setting $\mathbf P_{m_0}=\emptyset$. Then (1), (2), and (5) 
are clearly satisfied, while 
\begin{align*}
& \mathcal E(\mathbf P \setminus\mathbf P_{m_0}) = \mathcal E(\mathbf P) \le
2^{m_0n}\\ & \mathcal M(\mathbf P \setminus\mathbf P_{m_0}) = \mathcal M(\mathbf P)
\le 2^{2m_0n} \, , 
\end{align*}
hence (3) and (4) are also satisfied for $\mathbf P_{m_0}$. 

Suppose we have selected pairwise disjoint sets 
$\mathbf P_{m_0}$, 
$\mathbf P_{m_0-1} , \dots ,   \mathbf P_{m }$ for some $m<m_0$ such that 
(1)-(5)  are satisfied for all $j\in \{m_0, m_0-1, \dots , m\}$. We 
will construct a set of tiles $\mathbf P_{m-1}$ disjoint 
from all the sets already constructed such that 
(1)-(5)  are satisfied for all $j=m-1$. This procedure is given by  
decreasing induction.  We will need to consider the following four cases. 

\begin{case}   $\mathcal E\big(\mathbf P \setminus (\mathbf P_{m_0}  
\cup \dots \cup \mathbf P_{m} )\big) \le 2^{(m-1)n}$ and
$\mathcal M\big(\mathbf P  \setminus(\mathbf P_{m_0}
\cup \dots \cup \mathbf P_{m} )\big) \le 2^{2(m-1)n}$.
\end{case}

In this case set $\mathbf P_{m-1}=\emptyset$ and observe that (1)-(5) 
trivially hold. 

\begin{case}   $\mathcal E\big(\mathbf P  \setminus (\mathbf P_{m_0}  
\cup \dots \cup \mathbf P_{m} )\big) > 2^{(m-1)n}$ and   
$\mathcal M\big(\mathbf P \setminus (\mathbf P_{m_0}  
\cup \dots \cup \mathbf P_{m} )\big) \le 2^{2(m-1)n}$.  
\end{case}

Use Lemma \ref{lem2}   to find a subset 
$\mathbf P_{m-1}$ of $\mathbf P \setminus (\mathbf P_{m_0}  
\cup \dots \cup \mathbf P_{m} )$ such that 
\begin{equation}\lab{411}
  \mathcal E\big(\mathbf P \setminus (\mathbf P_{m_0}  
\cup \dots \cup \mathbf P_{m} \cup \mathbf P_{m-1})\big) \le 
\f12 \mathcal E\big(\mathbf P \setminus (\mathbf P_{m_0}  
\cup \dots \cup \mathbf P_{m} )\big)\le \f12 2^{mn}  
\end{equation}
and $\mathbf P_{m-1}$ is a union of trees (whose set of 
tops we   denote by $\mathbf P_{m-1}^*$) such that  
\begin{equation}\lab{421}
\sum_{t\in \mathbf P_{m-1}^*} |I_t| \le 
C_2 \mathcal E\big(\mathbf P \setminus (\mathbf P_{m_0}  
\cup \dots \cup \mathbf P_{m} )\big)^{-2} \le C_2 2^{-2(m-1)n}. 
\end{equation}
Then   (\ref{411}) gives (3)   and 
  (\ref{421}) gives (5) for $j=m\!-\! 1$. Since 
$$
\mathcal E(\mathbf P_{m-1}) \le \mathcal E\big(\mathbf P \setminus (\mathbf P_{m_0}  
\cup \dots \cup \mathbf P_{m} )\big) \le 2^{mn} = 2^{((m-1)+1)n}\, , 
$$
estimate (1) is satisfied for  $j=m\!-\! 1$. Also by our induction 
hypothesis we have 
$$
  \mathcal M\big(\mathbf P \setminus (\mathbf P_{m_0}  
\cup \dots \cup \mathbf P_{m} \cup \mathbf P_{m-1})\big) \le 
  \mathcal M\big(\mathbf P \setminus (\mathbf P_{m_0}  
\cup \dots \cup \mathbf P_{m} )\big)\le   2^{2(m-1)n} \, ,  
$$
hence (4) is satisfied for $j=m\!-\! 1$.  Finally 
$$
\mathbf P_{m-1} \subset \mathbf P \setminus (\mathbf P_{m_0}  
\cup \dots \cup \mathbf P_{m} )
$$
 and hence its mass is at most 
the mass of the latter which is trivially bounded by 
$2^{(2(m-1)+2)n}$, thus
(2) is  also satified  for $j=m\!-\! 1$.

\begin{case}    $\mathcal E\big(\mathbf P  \setminus (\mathbf P_{m_0}  
\cup \dots \cup \mathbf P_{m} )\big) \le 2^{(m-1)n}$ and   
$\mathcal M\big(\mathbf P \setminus (\mathbf P_{m_0}  
\cup \dots \cup \mathbf P_{m} )\big) > 2^{2(m-1)n}$.  
\end{case}

In this case we repeat the argument in case 2 with the roles of 
the mass and energy reversed. Precisely, use 
Lemma \ref{lem1}   to find a subset  
$\mathbf P_{m-1}$ of the set $\mathbf P \setminus (\mathbf P_{m_0}  
\cup \dots \cup \mathbf P_{m} )$ such that 
\begin{equation}\lab{431}
  \mathcal M\big(\mathbf P \setminus (\mathbf P_{m_0}  
\cup \dots \cup \mathbf P_{n} \cup \mathbf P_{m-1})\big) \le 
\f14 \mathcal M\big(\mathbf P \setminus (\mathbf P_{m_0}  
\cup \dots \cup \mathbf P_{m} )\big)\le \f14 2^{2mn}  
\end{equation}
and $\mathbf P_{m-1}$ is a union of trees (whose set of 
tops we   denote by $\mathbf P_{m-1}^*$) such that  
\begin{equation}\lab{441}
\sum_{t\in \mathbf P_{m-1}^*} |I_t| \le 
C_1 \mathcal M\big(\mathbf P \setminus (\mathbf P_{m_0}  
\cup \dots \cup \mathbf P_{m} )\big)^{-1} \le C_1 2^{-2(m-1)n}. 
\end{equation}
Then  (\ref{431}) gives (4)   and 
(\ref{441}) gives (5) for $j=m\!-\! 1$. By induction we have 
$$
\mathcal M(\mathbf P_{m-1}) \le \mathcal M\big(\mathbf P \setminus (\mathbf P_{m_0}  
\cup \dots \cup \mathbf P_{m} )\big) \le 2^{2mn}
= 2^{(2(m-1)+2)n}\, , 
$$
thus (2) is satisfied for $j=m\!-\! 1$.  Finally (1) and (3) 
follow from the inclusion $\mathbf P_{m-1}\subset 
\mathbf P \setminus (\mathbf P_{m_0}  
\cup \dots \cup \mathbf P_{m} )$ and the assumption 
$\mathcal E\big(\mathbf P  \setminus (\mathbf P_{m_0}  
\cup \dots \cup \mathbf P_{m} )\big) \le 2^{(m-1)n}$. 
This concludes the proof of (1)-(5) for $j=m\!-\! 1$.

\begin{case}  $\mathcal E\big(\mathbf P  \setminus (\mathbf P_{m_0}  
\cup \dots \cup \mathbf P_{m} )\big) > 2^{(m-1)n}$ and   
$\mathcal M\big(\mathbf P \setminus (\mathbf P_{m_0}  
\cup \dots \cup \mathbf P_{m} )\big) > 2^{2(m-1)n}$.  
\end{case}

This is the most difficult case since it involves elements 
from both of the previous cases. 
We start by using 
Lemma \ref{lem1}   to find a subset  
$\mathbf P_{m-1}'$ of the set $\mathbf P \setminus (\mathbf P_{m_0}  
\cup \dots \cup \mathbf P_{m} )$ such that 
\begin{equation}\lab{451}
  \mathcal M\big(\mathbf P \setminus (\mathbf P_{m_0}  
\cup \dots \cup \mathbf P_{m} \cup \mathbf P_{m-1}')\big) \le 
\f14 \mathcal M\big(\mathbf P \setminus (\mathbf P_{m_0}  
\cup \dots \cup \mathbf P_{m} )\big)\le \f14 2^{2mn}  
\end{equation}
and $\mathbf P_{m-1}'$ is a union of trees (whose set of 
tops we   denote by $(\mathbf P_{m-1}')^*$) such that  
\begin{equation}\lab{461}
\sum_{t\in (\mathbf  {P}_{m-1}')^*} |I_t| \le 
C_1 \mathcal M\big(\mathbf P \setminus (\mathbf P_{m_0}  
\cup \dots \cup \mathbf P_{m} )\big)^{-1} \le C_1 2^{-2(m-1)n}. 
\end{equation}
We now consider the following two subcases of case 4.

\begin{subcase4a}  
$\mathcal E\big(\mathbf P  \setminus (\mathbf P_{m_0}  
\cup \dots \cup \mathbf P_{m} \cup \mathbf P_{m-1}')\big) \le 2^{(m-1)n}$ 
\end{subcase4a}

In this subcase, we set $\mathbf P_{m-1}=\mathbf P_{m-1}'$. Then 
(3) is automatically satisfied for $j=m\!-\! 1$ and also (5) is 
satisfied in view of (\ref{461}). By the inductive hypothesis 
we have 
$\mathcal E\big(\mathbf P  \setminus (\mathbf P_{m_0}  
\cup \dots \cup \mathbf P_{m}  )\big) \le 2^{mn }=2^{((m-1)+1)n}$  
and also 
$\mathcal M\big(\mathbf P  \setminus (\mathbf P_{m_0}  
\cup \dots \cup \mathbf P_{m}  )\big) \le 2^{2mn}=2^{(2(m-1)+2)n}$. 
Since $\mathbf P_{m-1}$ is contained in 
$\mathbf P  \setminus (\mathbf P_{m_0}  
\cup \dots \cup \mathbf P_{m}  )$ the same estimates hold for 
$\mathcal E(\mathbf P_{m-1})$ and $\mathcal M(\mathbf P_{m-1})$, thus (1) and (2) 
also hold for $j=m-1$.  Finally (4) for $j=m-1$ follows from 
(\ref{45}) since $\mathbf P_{m-1}'=\mathbf P_{m-1}$. 

\begin{subcase4b}  
$\mathcal E\big(\mathbf P  \setminus (\mathbf P_{m_0}  
\cup \dots \cup \mathbf P_{m} \cup \mathbf P_{m-1}')\big) > 2^{(m-1)n}$ 
\end{subcase4b}

Here we use Lemma \ref{lem2} one more time  to find a subset 
$\mathbf P_{m-1}''$ of the set $\mathbf P \setminus (\mathbf P_{m_0}  
\cup \dots \cup \mathbf P_{m}  \cup \mathbf P_{m-1}')$ such that 
\begin{align}\begin{split}\lab{471}
 & \mathcal E\big(\mathbf P \setminus (\mathbf P_{m_0}  
\cup \dots \cup \mathbf P_{m} \cup \mathbf P_{m-1}'
\cup \mathbf P_{m-1}'')\big) \\ \le  &
\f12 \mathcal E\big(\mathbf P \setminus (\mathbf P_{m_0}  
\cup \dots \cup \mathbf P_{m}\cup \mathbf P_{m-1}' )\big)    
\end{split}\end{align}
and $\mathbf P_{m-1}''$ is a union of trees (whose set of 
tops we   denote by $(\mathbf {P}_{m-1}'')^*$) such that  
\begin{equation}\lab{481}
\sum_{t\in (\mathbf  {P}_{m-1}'')^*} |I_t| \le 
C_2 \mathcal E\big(\mathbf P \setminus (\mathbf P_{m_0}  
\cup \dots \cup \mathbf P_{m}
\cup \mathbf P_{m-1}' )\big)^{-2} \le C_2 2^{-2(m-1)n}. 
\end{equation}
We set $\mathbf P_{m-1} = \mathbf P_{m-1}' \cup \mathbf P_{m-1}''$
and we observe that $\mathbf P_{m-1}$ is disjoint from all the 
previously selected $\mathbf P_{j}$'s. 
Since by the  induction hypothesis  the last term in (\ref{471}) 
is bounded by 
$\f12 \mathcal E\big(\mathbf P \setminus (\mathbf P_{m_0}  
\cup \dots \cup \mathbf P_{m}  )\big) \le \f12 2^{mn}$,  
the first term in (\ref{471}) 
is also bounded by $2^{(m-1)n}$, thus (3) holds for $j= m\! - \! 1$. 
Likewise, since 
\begin{align*}
&\mathcal E(\mathbf P_{m-1} ) \le \mathcal E\big(\mathbf P \setminus (\mathbf
P_{m_0}  
\cup \dots \cup \mathbf P_{m}  )\big) \le 2^{mn} = 2^{((m-1)+1)n} \\ 
&\mathcal M(\mathbf P_{m-1} ) \le \mathcal M\big(\mathbf P \setminus (\mathbf P_{m_0}  
\cup \dots \cup \mathbf P_{m}  )\big) \le 2^{2mn} = 2^{(2(m-1)+2)n}\, , 
\end{align*}
(1) and (2) are satified for  $j= m\! - \! 1$. Since 
$$
\mathcal M\big(\mathbf P \setminus (\mathbf P_{m_0}  
\cup \dots \cup \mathbf P_{m }\cup \mathbf P_{m-1}  )\big)\le 
\mathcal M\big(\mathbf P \setminus (\mathbf P_{m_0}  
\cup \dots \cup \mathbf P_{m }\cup \mathbf P_{m-1}'  )\big)\, 
$$
(\ref{451}) implies that (4) is satisfied for $j= m\! - \! 1$. 
Now each of $\mathbf P_{m-1}'$ and $ \mathbf P_{m-1}''$ is 
given as a union of trees, thus the same is true for $\mathbf P_{m-1}$. 
The set of tops of all of these trees, call it $(\mathbf P_{m-1})^* $,
is contained in the union of the set of tops of the trees in  
$\mathbf P_{m-1}'$ and the trees in $ \mathbf P_{m-1}''$, i.e in 
$(\mathbf {P}_{m-1}')^*\cup (\mathbf {P}_{m-1}'')^*$. This 
implies that 
\begin{align*}
\sum_{t\in (\mathbf {P}_{m-1})^*} |I_t|  \le &
\sum_{t\in (\mathbf {P}_{m-1}')^*} |I_t|  +
\sum_{t\in (\mathbf {P}_{m-1}'')^*} |I_t| \\
\le & (C_1+C_2) 2^{-2(m-1)n}
=C_0 2^{-2(m-1)n}
\end{align*}
in view of (\ref{461}) and 
(\ref{481}). 
This proves (5) for $j= m\! - \! 1$ and concludes the inductive 
step   $j= m\! - \! 1$. The 
construction of the $\mathbf P_j$'s is now complete.

\section{Proof of lemma 1}
\setcounter{equation}0 

Given a finite set of tiles $\pp$, set $\mu = \cm (\pp)$ and define 
$$
\pp'=\{p \in \pp: \cm(\{p\}) > \f {1} {4} \mu\}.
$$
Clearly $\cm(\pp \setminus \pp') \le \f {1} {4} \mu$, thus it remains to show that 
$\pp'$ satisfies
(\ref{42}).  By definition of the mass, for each $p \in \pp'$ there is a tile
 $u(p)=u \in \dd$ such that
\begin{equation}\lab{51}
 \intl_{E \cap \N^{-1}[\om_u]} \df {|\iuu|^{-1}} 
{\Bigl(1+ \f {|\x-c(\iuu)|} {|\iuu|^{1/n}}\Bigr)^{10n}} \, d\x > \f {\mu} {4}.
\end{equation}
Set $\bu=\{u(p):p \in \pp'\}$, and let $\bu_{\textup {max}}$ be the subset of $\bu$
containing all
maximal elements of $\bu$ under the partial order on tiles.  As observed earlier,
the tiles in $\bu$ can be grouped into trees with
tops in  $\bu_{\textup {max}}$.  Now $\bu$ is not necessarily a subset of $\pp'$, but each
$u \in \bu$ is associated to a $p \in \pp'$ as described above.  In particular, if $p$ is a
maximal element in $\pp'$, then there exists a $u \in \bu$ with $p<u$ such that
(\ref{51}) holds.  If this $u$ is not in $\bu_{\textup{max}}$ then there exists
$u' \in \bu_{\textup {max}}$
with $u < u'$. We must then have $u'$ associated to another $p' \in \pp'$ which implies, by
maximality of $p$, that
$p' < p$.  Hence for each maximal element  $p \in \pp'$ there exists a unique element
$u \in \bu_{\textup {max}}$ with $p < u$, and there is at most one such maximal element for each
 $u \in \bu_{\textup {max}}$.  Therefore, we will show
\begin{equation}\label{52}
\summ_{u \in \bu_{\textup {max}}} |\iuu| \le C_1 \mu^{-1},
\end{equation}
which implies (\ref{42}).
Now we will rewrite (\ref{51}) as
$$
\f {2^n-1}{2^{n+2}} \mu \summ_{k=0}^\nf 2^{-kn} <  \summ_{k=0}^\nf \ \intl_{\substack {E \cap \N^{-1}[\om_u] \\ \cap 
(2^k \iuu \setminus 2^{k-1} \iuu)}} \df {|\iuu|^{-1}}
{\Bigl(1+ \f {|\x-c(\iuu)|} {|\iuu|^{1/n}}\Bigr)^{10n}} \, d\x,
$$
where we set $\f {1} {2} \iuu = \emptyset$.  This estimate holds for all $u \in \bu$, so in particular for 
every $u \in \bu_{\textup {max}}$ there exists a $k \ge 0$ such that                    
\begin{align*}
\f {2^n-1}{2^{n+2}} \mu  |\iuu| 2^{-kn} < & \  \intl_{\substack {E \cap \N^{-1}[\om_u] \\ \cap
(2^k \iuu \setminus 2^{k-1} \iuu)}} \df {1} 
{\Bigl(1+ \f {|\x-c(\iuu)|} {|\iuu|^{1/n}}\Bigr)^{10n}} \, d\x \\
\le & \df {|E \cap \N^{-1}[\om_u] \cap
(2^k \iuu \setminus 2^{k-1} \iuu)|} {C(\sqrt {n})^{10n} 2^{10kn}},
\end{align*}
where the second inequality above follows from the fact that 
$$|\iuu|^{-1/n}|\x-c(\iuu)|  \sim \sqrt {n}2^k$$ 
for $\x \in (2^k \iuu \setminus 2^{k-1} \iuu)$.  Here and throughout the paper $C$
denotes a constant depending only on dimension and whose value
may change at different places in the proof.
Now we define for $k \ge 0$
$$
\bu_k=\{ u \in \bu_{\textup {max}}:C \mu |\iuu| 2^{9kn} < |E \cap \N^{-1}[\om_u]
\cap 2^k \iuu | \}.
$$
Since  $\bu_{\textup {max}} = \bigcup _{k=0}^\nf \bu_k$, if we show that
\begin{equation}\label{0}
\summ_{u \in \bu_k} |\iuu| \le C 2^{-8kn} \mu^{-1}, 
\end{equation}
summing over $k \in \zzz^+$ gives us estimate (\ref{52}).

We now concentrate on showing estimate (\ref{0}). Fix $k \ge 0$ and select an element $v_0 \in \bu_k$ 
so that 
$|I_{v_0}|$ is largest possible.  Then select an element 
$v_1 \in \bu_k \setminus\{v_0\}$ such that 
the enlarged rectangle $(2^k I_{v_1}) \times \om_{v_1}$ is disjoint from  
$(2^k I_{v_0}) \times \om_{v_0}$ and $|I_{v_1}|$ is largest possible. 
 Continuing by induction, at 
the $j$-th step we select an element $v_j \in \bu_k \setminus\{v_0, \ldots, v_{j-1}\}$
 so that
$(2^k I_{v_j}) \times \om_{v_j}$ is disjoint from the enlarged rectangles 
of previously selected tiles
and $|I_{v_j}|$ is largest possible.  Since we have a finite set of tiles,
 this process will 
terminate, and we will have the set of selected tiles in $\bu_k$, 
which we will call $\bv_k$.

Next we make some key observations about the tiles.  First, note that elements of 
$\bu_k$ are maximal
in $\bu$ and therefore disjoint.  Second, for any $u \in \bu_k$ there exits a selected
 tile $v \in 
\bv_k$ with $|\iuu| \le |I_v|$ and such that the enlarged rectangles of $u$ and $v$ intersect.  We will 
associate $u$ to this $v$.  Third, if $u$ and $u'$ are both associated to the 
same $v$, then $\iuu$
and $I_{u'}$ are disjoint.  Indeed, $(2^k I_u) \times \om_u$ intersects 
$(2^k I_v) \times \om_v$
which means $2^k I_u \cap 2^k I_v \neq \emptyset$ and 
$\om_u \cap \om_v \neq \emptyset$. This implies, 
together
with the fact that  $|\iuu| \le |I_v|$, that $\om_u \supset \om_v$.  
Similarily $\om_{u'} \supset 
\om_v$.  Therefore, one of $\om_u$ and $\om_{u'}$ contains the other.
But $u$ and $u'$ are disjoint,
thus $\iuu$ is disjoint from $I_{u'}$.  Finally , all tiles 
$u \in \bu_k$ associated to a particular
$v \in \bv_k$ satisfy $\iuu \subset 2^{k+2} I_v$.

From the observations above and the definition of $\bu_k$, we have
{\allowdisplaybreaks
\begin{align*}
\summ_{u \in \bu_k} |\iuu| \le & \summ _{v \in \bv_k} 
\summ _{\substack {u \in \bu_k \\ u \,
\textup{assoc}\, v}} |\iuu| \\ 
= & \summ _{v \in \bv_k} \biggl| \bigcup_{\substack
{u \in \bu_k \\ u \, \textup{assoc}\, v}} 
\iuu \biggr|  \le  \summ _{v \in \bv_k}2^{(k+2)n} |I_v| \\
\le & C \mu^{-1}2^{-9kn} 2^{(k+2)n} 
\summ _{v \in \bv_k} |E \cap \N^{-1}[\om_v] \cap 2^k I_v | \\
 \le & C 2^{2n} \mu^{-1}2^{-8kn}|E|  \le  C 2^{2n} \mu^{-1}2^{-8kn},
\end{align*}}
where we have used that for $v \in \bv_k$, the enlarged rectangles are disjoint, 
and therefore so are 
the subsets $E \cap \N^{-1}[\om_v] \cap 2^k I_v$ of $E$.

\section{Proof of lemma 2}
\setcounter{equation}0 

We begin by fixing a finite set of tiles $\pp$ and  $r = 2^n$.
This choice of $r$ ensures that ${\pmb {\eta}} < {\pmb {\xi}}$ in the lexicographical order
for all ${\pmb {\eta}} \in \om_{p(1)}$ and
${\pmb {\xi}} \in \om_{p(r)}$.
For $r \ne 2^n$ the proof goes
through by a suitable permutation of the coordinates of $\mathbb R^{n}$
which changes the coordinate that takes precedence in the lexicographical order.
Here we note that we can be less precise by taking any linear functional $L$ that
separates $\om_{p(1)}$ and $\om_{p(r)}$ in any given cube $\om_p$.  Then we let
${\pmb {\eta}} < {\pmb {\xi}}$ if $L({\pmb {\eta}}) < L({\pmb {\xi}})$.  In particular,
we can take $L$ to be the projection onto the appropriate axis so that the usual linear ordering
on $\mathbb R$ is relevant.
Let $\epsilon$ denote $\ce(\pp)$. Define for $\tr'$ a $2^n$-tree
$$
\De (\tr')= \biggl(|\itr|^{-1} \summ_{p \in \tr'} |\lng f, \phi_p \rng|^2 \biggr)^{1/2}.
$$
Consider all $2^n$-trees $\tr'$ contained in $\pp$ which satisfy
\begin{equation}\lab{61}
\De (\tr') \ge \f {1} {2} \epsilon.
\end{equation}
Among these select a $2^n$-tree $\tr_1'$ such that
$c(\om_{\tr_1'})$ is minimal in the lexicographical order.
  Let $\tr_1$ be the set of all $p \in \pp$ such that 
$p<p_{\tr_1'}=p_{\tr_1}$.  
 In other words $\tr_1$ is the maximal tree containing $\tr'_1$ with the same top as 
 $\tr'_1$.
 Now consider all $2^n$-trees contained in $\pp \setminus \tr_1$ 
and select a
$2^n$-tree $\tr_2,$ such that $c(\om_{\tr_2'})$ is minimal. 
Let $\tr_2 $ be the set of all $p \in \pp$ such that $p<p_{\tr_2'}=p_{\tr_2}$. 
Continue inductively to obtain a finite sequence of pairwise disjoint $2^n$-trees
$
\tr_1' \, ,\, \tr_2'\, ,\, \ldots\, ,\, \tr_q'
$
and pairwise disjoint trees
$
\tr_1\, ,\, \tr_2\, ,\, \ldots\, ,\, \tr_q
$
where $p_{\tr_j'}=p_{\tr_j}$, $\tr_j' \subset \tr_j$,
 and the $\tr_j'$ satisfy (\ref{61}).  Let 
$$
\pp''= \bigcup_{j=1}^q  \tr_j \ ,
$$
then clearly 
$$
\ce (\pp \setminus \pp'') \le \f {1} {2} \epsilon \ .
$$ 
Thus we need to show that $\pp''$ satisfies condition (\ref{44}) of Lemma \ref{lem2}, i.e.
\begin{equation}\lab{64}
\summ_{j=1}^q  |I_{\tr_j}| \le \f {C_2} {\epsilon^2}.
\end{equation}

Since the trees $\tr_j'$ satisfy (\ref{61}) and $\|f\|_{L^2(\rn)}=1$,
{\allowdisplaybreaks
\begin{align}\begin{split}\label{62}
  \f {1} {4} \epsilon^2 \summ_j |I_{\tr_j}|  
\le & \summ_j \summ_{p \in \tr_j'} | \lng f,\phi_p \rng|^2 \\ 
 = &  \summ_j \summ_{p \in \tr_j'} \lng f,\phi_p \rng \overline{\lng f,\phi_p \rng} \\
= & \lng f,\summ_j \summ_{p \in \tr_j'} \overline{\lng f, \phi_p \rng}\phi_p \rng \\
\le &\biggl\| \summ_j \summ_{p \in \tr_j'} \lng f, \phi_p \rng \phi_p \biggr\|_{L^2(\rn)}.\\
\end{split}\end{align}}
Letting $\bu =\bigcup_j \tr_j'$, we will show that
\begin{equation}\lab{63}
\biggl\| \summ_{p \in \bu} \lng f, \phi_p \rng \phi_p \biggr\|_{L^2(\rn)} \le C
\Bigl( \epsilon^2 \summ_j 
|I_{\tr_j}| \Bigr)^{1/2},
\end{equation}
which, together with (\ref{62}), will give us (\ref{64}). The square of the left hand side of 
(\ref{63}) can be estimated by
\begin{equation}\lab{65}
\summ_{\substack{ p,u \in \bu \\ \oms=\omu }} | \lng f, \phi_p \rng \lng f,\phi_u \rng 
\lng \phi_p,\phi_u \rng | + 2 \summ_{\substack{ p,u \in \bu \\ 
\oms \subset \om_{u(1)} }}
| \lng f, \phi_p\rng \lng f,\phi_u \rng  \lng \phi_p,\phi_u \rng |.
\end{equation}
Here we have used that $\lng \phi_p,\phi_u \rng = 0$ unless $\omsone$ intersects $\om_{u(1)}$
 which implies
that either $\oms=\omu$ or $\om_{u(1)}$ contains $\oms$ or $\omsone$
contains $\omu$.
We are then able
 to utilize the symmetry in $p$ and $u$ to combine the off diagonal terms.
We estimate $|\lng f, \phi_p \rng|$ and $|\lng f,\phi_u \rng|$
by the larger one and bound the first term in (\ref{65}) by
{\allowdisplaybreaks 
\begin{align*}
& \summ_{p \in \bu} | \lng f, \phi_p \rng |^2 \summ_{\substack{ u \in \bu \\ 
\oms=\omu }}
|\lng \phi_p,\phi_u \rng| \\
\le & \summ_{p \in \bu} | \lng f, \phi_p \rng |^2 \summ_{\substack{ u \in \bu \\ 
\oms=\omu }}
C \df {\min \left( \f {|\iuu|} {|\is|},\f {|\is|} {|\iuu|} \right)^{1/2}}
{\left( 1 + \f {|c(\is) - c(\iuu)|} {\max (|\is|,|\iuu|)^{1/n}} \right)^{10n}}  \\
\leq & \summ_{p \in \bu} | \lng f, \phi_s \rng |^2 \summ_{\substack{ u \in \bu \\ 
\oms=\omu }}
C \int_{\iuu} \f {1} {|\is|} \biggl( 1 + \f {|\x-c(\is)|} {|\is|^{1/n}} \biggr)^{-10n}
 \ d\x \\
\leq & C \summ_{p \in \bu} | \lng f, \phi_p \rng |^2 \\
= & C \sum_j \sum_{p \in \tr_j'}|\itj| |\itj|^{-1} | \lng f, \phi_p \rng |^2 \\
\leq & C \sum_j |\itj|  \epsilon^2 ,
\end{align*}}
where we have used that for $p \in \bu$, the $\iuu$ for which $\oms=\omu$ 
are pairwise disjoint.

Using Cauchy-Schwarz, the second term in (\ref{65}) can be estimated by 
{\allowdisplaybreaks
\begin{align*}
& 2 \sum_j \sum_{p \in \tr_j'} |\lng f,\phi_p \rng|
\sum_{\substack{ u \in \bu \\ \oms \subset \om_{u(1)} }}
|\lng f,\phi_u \rng|  | \lng \phi_p,\phi_u \rng | \\
\leq & 2 \sum_j \left\{ \sum_{p \in \tr_j'}|\lng f,\phi_p \rng|^2 \right\}^{1/2}
\left\{\sum_{p \in \tr_j'} \Biggl( \, \sum_{\substack{ u \in \bu \\
 \oms \subset \om_{u(1)} }}
|\lng f,\phi_u \rng|  | \lng \phi_p,\phi_u \rng |\biggr)^2 \right\}^{1/2}\\
\leq & 2 \sum_j |\itj|^{1/2} \De(\tr_j') \left\{\sum_{p \in \tr_j'} \Biggl(
\, \sum_{\substack{ u \in \bu \\ \oms \subset \om_{u(1)} }}
|\lng f,\phi_u \rng|  | \lng \phi_o,\phi_u \rng | \Biggr)^2 \right\}^{1/2}\\
\leq & 2 \epsilon \sum_j |\itj|^{1/2} \left\{\sum_{p \in \tr_j'}
\Biggl( \, \sum_{\substack{ u \in \bu \\ \oms \subset \om_{u(1)} }}
|\lng f,\phi_u \rng|  | \lng \phi_p,\phi_u \rng | \Biggr)^2 \right\}^{1/2}.\\
\end{align*}}
To complete the proof, we need to show that the expression inside the 
curly brackets is bounded by
$C \epsilon^2 |\itj|$.  Since for a single tile $u$
$$
\ce(\{u\}) = \left( |\iuu|^{-1} |\lng f,\phi_u \rng|^2 \right)^{1/2} =
|\iuu|^{-1/2} |\lng f,\phi_u \rng| \leq \epsilon,
$$  
we get that
$$
\sum_{p \in \tr_j'} 
\Biggl( \, \sum_{\substack{ u \in \bu \\ \oms \subset \Om_{u(1)} }}
|\lng f,\phi_u \rng|  | \lng \phi_p,\phi_u \rng | \Biggr)^2
\leq \epsilon^2 \sum_{p \in \tr_j'} 
\Biggl( \, \sum_{\substack{ u \in \bu \\ \oms \subset \om_{u(1)} }}
|\iuu|^{1/2} | \lng \phi_p,\phi_u \rng | \Biggr)^2.
$$
Thus we now need to show that
\begin{equation}\lab{67}
\sum_{p \in \tr_j^r} 
\Biggl( \, \sum_{\substack{ u \in \bu \\ \oms \subset \om_{u(1)} }}
|\iuu|^{1/2} | \lng \phi_p,\phi_u \rng | \Biggr)^2
\leq C |\itj|.
\end{equation}
To prove this, we will need the following lemma. 

\noindent \begin{lemma}\lab{lem4}
Let $p \in \tr_j'$ and $u \in \tr_k'$.  Then if $\oms \subset \om_{u(1)}$,
 we have $\iuu \cap  \itj
= \emptyset$.  If $u \in \tr_k' , v \in \tr_l', u \neq v, \oms \subset \om_{u(1)}$, 
and
$\oms \subset \om_{v(1)}$ for some fixed $p \in \tr_j'$, 
then $ \iuu \cap I_v = \emptyset$.  
\end{lemma}

\begin{proof}
  Since $\oms \subset \om_{u(1)}$ and $\tr_j'$ and $\tr_k'$ are $2^n$-trees, 
$\tr_j'$ and $\tr_k'$ are not the same tree.  Otherwise $\oms \subset \om_{u(2^n)}$.
We know that $\om_{\tr_j'} \subset \oms \subset \om_{u(1)}$, which implies that
$c(\om_{\tr_j'})$ is contained in $\om_{u(1)}$. We also have
$\om_{\tr_k'} \subset \om_{u(2^n)}$, which implies that
$c(\om_{\tr_k'})$is contained in $\om_{u(2^n)}$.
Therefore $c(\om_{\tr_j'}) \leq c(\om_{\tr_k'})$
in the lexicographical order 
which means that $\tr_j'$ 
was chosen before $\tr_k'$ in the original selection process.  
Now suppose $\iuu \cap  \itj
\neq \emptyset$. Then either $\iuu \subset  \itj$ or $\iuu \supset  \itj$, 
however $\om_{\tr_j} 
\subset \omu$ implies that $\iuu \subset  \itj$.   Thus we have  
$\om_{\tr_j} 
\subset \omu$ and $\iuu \subset I_{\tr_j}$ which says that
$u$ belongs to the tree $\tr_j$.  However, $u \in \tr_k$ and thus was chosen from 
$\pp \setminus 
\tr_j$, which gives a contradiction.  Thus 
$\iuu \cap \itj= \emptyset$.

Next suppose that $u \in \tr_k' , v \in \tr_l', u \neq v$, and
$\oms \subset (\om_{u(1)} \cap \om_{v(1)})$ for some fixed $p \in \tr_j'$.
We have three cases to consider: (a) $\omu \subset \om_{v(1)}$ which means 
$I_v \cap I_{\tr_k}= \emptyset$ and thus $I_v \cap \iuu =\emptyset$, 
(b) $\om_v \subset \om_{u(1)}$
 which means $I_u \cap  I_{\tr_l}= \emptyset$ and thus 
 $I_v \cap \iuu = \emptyset$, and (c) 
$\om_v = \om_u$
which tells us $|\iuu|=|I_v|$, thus $I_u$ and $I_v$ are disjoint since 
$u$ and $v$ don't coincide. 
\end{proof}

We now return to estimate (\ref{67}).  Observe that Lemma \ref{lem4} 
tells us that for the
tiles $u \in \bu$ appearing in the interior sum of (\ref{67}), the 
$\iuu$ are pairwise disjoint and 
contained in $(I_{\tr_j})^c$.  Thus we have
{\allowdisplaybreaks
\begin{align*}
& \sum_{p \in \tr_j'} 
\Biggl( \, \sum_{\substack{ u \in \bu \\ \oms \subset \om_{u(1)} }}
|\iuu|^{1/2} | \lng \phi_p,\phi_u \rng | \Biggr)^2 \\
\le & C \sum_{p \in \tr_j'}
\Biggl( \, \sum_{\substack{ u \in \bu \\ \oms \subset \om_{u(1)} }}
|\iuu|^{1/2} \biggl( \f {|\is|} {|\iuu|} \biggr)^{1/2}
\int_{\iuu} \f {|\is|^{-1}} {\Bigl( 1 + \f {|\x-c(\is)|} 
{|\is|^{1/n}} \Bigr)^{10n}}  \, d\x
\Biggr)^2 \\
\le & C \sum_{p \in \tr_j'}|\is| 
\Biggl( \, \sum_{\substack{u \in \bu \\ \oms \subset \om_{u(1)} }}
\int_{\iuu} \f {|\is|^{-1}} {\Bigl( 1 + \f {|\x-c(\is)|} 
{|\is|^{1/n}} \Bigr)^{10n}}  \, d\x
\Biggr)^2 \\
\le & C \sum_{p \in \tr_j'}|\is| \Biggl( \,
\int_{(\itj)^c} \f {|\is|^{-1}} {\Bigl( 1 + \f {|\x-c(\is)|} 
{|\is|^{1/n}} \Bigr)^{10n}}  \, d\x
\Biggr)^2 \\
\le & C \sum_{p \in \tr_j'}|\is| 
\int_{(\itj)^c} \f {|\is|^{-1}} {\Bigl( 1 + \f {|\x-c(\is)|} 
{|\is|^{1/n}} \Bigr)^{10n}}  \, d\x \\
\le & C \sum_{k=0}^\nf 2^{-kn} |\itj| 
\sum_{\substack{ p \in \tr_j' \\ |\is|=2^{-kn} 
|\itj|}}    
\int_{(\itj)^c} \f {|\is| ^{-1}} 
{\Bigl( 1 + \f {|\x-c(\is)|} {|I_p|^{1/n}} \Bigr)^{10n}}  
\, d\x .\\
\end{align*}}
The proof of Lemma \ref{lem2} will be complete if we can show that
$$
\sum_{\substack{ p \in \tr_j' \\ |\is|=2^{-kn} 
|\itj|}}    
\int_{(\itj)^c} \f {|\is| ^{-1}} {\Bigl( 1 + \f {|\x-c(\is)|}
 {|I_{p}|^{1/n}} \Bigr)^{10n}}
\, d\x  \lesssim 2^{k(n-1)},
$$
thus allowing the sum in $k$ to converge.  Throughout the paper, $A \lesssim B$ means
that $A$ is less than or equal to $B$ up to a constant depending only on dimension.
The first observation we have is that 
$$
\int_{(\itj)^c} \f {|\is| ^{-1}} {\Bigl( 1 + \f {|\x-c(\is)|} 
{|I_{p}|^{1/n}} \Bigr)^{10n}}  
\, dx  \lesssim \left( \f {\textup{dist}((\itj)^c,I_p)} { |I_p|^{1/n}} \right)^{-9n}.
$$
To see this, note that by a change of variables, it suffices to let the center
 of $\itj$ be at the
origin.  Also note that we have the inequality
$$
\Bigl( 1 + \f {|\x-c(\is)|}{|I_{p}|^{1/n}}\Bigr)^{10n} \ge \prod_{i=1}^n
{\Bigl( 1 + \f {|x_i-c(\is)_i|}
{|I_{p}|^{1/n}} \Bigr)^{10}}.
$$
Therefore, the integral above is bounded by a constant times
\begin{align*}
 & \prod_{i=1}^n  \left(
\intl_{|x_i| > \f {1} {2} |\itj|}
\f {|\is| ^{-1/n}} {\Bigl( 1 + \f {|x_i-c(\is)_i|}
{|I_{p}|^{1/n}} \Bigr)^{10}} \, dx_i \right)\\
 \lesssim &\prod_{i=1}^n \left( \f {\textup{dist}((\itj)^c,I_p)} { |I_p|^{1/n}} \right)^{-9}
 \lesssim  \left( \f {\textup{dist}((\itj)^c,I_p)} { |I_p|^{1/n}} \right)^{-9n},
\end{align*}
where we have used that $|x_i-c(I_p)_i| \ge \textup{dist}((\itj)^c,I_p)$ for all
$i=1,2, \ldots, n$.
Now we need to sum over $p$ for a fixed scale $k$.  Consider an $n-1$
dimensional face of
$\itj$  and fix a cube $I_p$ whose face is contained in the face of $\itj$.
We allow the remaining coordinate to vary and sum over those $I_p$ in
this ``column''.
  In a fixed column, the
distances from $(\itj)^c$ to each
$I_p$ sum as
additive multiples of $|I_p|^{1/n}$.  For each face, there are $2^{k(n-1)}$
such columns.  Thus
\begin{align*}
 \sum_{\substack{ p \in \tr_j' \\ |\is|=2^{-kn}|\itj|}}
\left( \f {\textup{dist}((\itj)^c,I_p)} { |I_p|^{1/n}} \right)^{-9n} 
\lesssim & 2^{k(n-1)} \times (\# \textup{of faces}) \sum_{m=0}^{\nf} \f {1} {m^{9n}} \\
\lesssim & 2^{k(n-1)}.
\end{align*}

\section{Proof of Lemma 3 - The Tree Inequality}
\setcounter{equation}0

Let $\mathcal J$ be the collection of all maximal dyadic cubes $J$ such that
$3J$ does not contain any $I_{p}$ with $ p \in {\bf{T}}$. Then $\mathcal J$ is
a partition of $\mathbb R^{n}$. 

We can write the left hand side of (\ref{45}) as follows, where
the terms $\alpha_{p}$ are phase factors of modulus 1 which make up
for the absolute value signs in (\ref{45}): 
\begin{displaymath}
|| \sum_{p \in {\bf{T}}} \alpha_{p} \langle f, \phi_{p} \rangle \psi_{p}^{{\bf{N}}}
   1_{E_{2p}} ||_{1} \leq \mathcal K_{1} + \mathcal K_{2}
\end{displaymath}
where 
\begin{align}
E_{2p} &:= E \cap {\bf {N}}^{-1}[\omega_{p(2^{n})}],\\
\mathcal K_{1} &:= \sum_{J \in \mathcal J} \sum_{p \in {\bf{T}}\,;\, |I_{p}|
\leq 2^{n}|J|} ||\langle f, \phi_{p} \rangle \psi_{p}^{{\bf{N}}} 1_{E_{2p}}
||_{L^{1}\left(J\right)}, \label{k1} \\
\mathcal K_{2} &:= \sum_{J \in \mathcal J} \left|\left|\sum_{p \in
{\bf{T}}\,;\, |I_{p}| > 2^{n}|J|}
\alpha_{p} \langle f, \phi_{p} \rangle \psi_{p}^{{\bf{N}}}
   1_{E_{2p}} \right|\right|_{L^{1}\left(J\right)}. \label{k2}
\end{align}
Let 
\begin{displaymath}
\epsilon = \mathcal E({\bf{T}}) \quad {\mbox{ and }} \quad \mu = \mathcal M({\bf{T}}).
\end{displaymath}

We begin with $\mathcal K_{1}$. For every $p \in {\bf{T}}$, $\{p\}$ is a $2^n$-tree contained 
in ${\bf{T}}$, and therefore
\[ \left| \langle f, \phi_{p}\rangle \right| \leq \epsilon |I_{p}|^{\frac{1}{2}}. \]

\allowdisplaybreaks{ 
\begin{align*}
\mathcal K_{1} &\leq C \epsilon \sum_{J \in \mathcal J} \sum_{\begin{subarray}{c}p \in {\bf{T}} \\ 
|I_{p}| \leq 2^{n} |J| \end{subarray}} |I_{p}| \intl_{J \cap E \cap \N^{-1}[\om_p]} 
\frac{|I_{p}|^{-1}}{\left( 1 + \frac{\left|{\bf{x}} - c\left(I_{p}\right)\right|}
{\left| I_{p}\right|^{\frac{1}{n}}}\right)^{20n}}\, d{\bf{x}}\\ 
& \leq C \epsilon \mu \sum_{J \in \mathcal J} \sum_{\begin{subarray}{c}p \in {\bf{T}} \\ 
|I_{p}| \leq 2^{n} |J| \end{subarray}} |I_{p}| \,  \sup_{{\bf{x}} \in J} 
\left( 1 + \frac{\left|{\bf{x}} - c\left(I_{p}\right)\right|}
{\left| I_{p}\right|^{\frac{1}{n}}}\right)^{\! \! -10n}\\
& \leq C \epsilon \mu \sum_{J \in \mathcal J} \sum_{k\, : \,2^{kn} \leq 2^{n} |J|} 2^{kn} 
\sum_{\begin{subarray}{c}p \in {\bf{T}} \\ |I_{p}|= 2^{kn}\end{subarray}}
\left( 1 + \frac{{\mbox{dist}}\left(J, I_{p}\right)}{2^{k}} \right)^{\! \! -10n},
\end{align*}}
where we have used that for ${\bf{x}} \in J$,
\[ |{\bf{x}} - c(I_{p})| \geq {\mbox{dist}}(J, c(I_{p})) \geq {\mbox{dist}}(J, I_{p}) + \frac{2^{k}}{2},\]
hence
\[ \frac{|{\bf{x}} - c(I_{p})|}{2^{k}} \geq \frac{{\mbox{dist}}(J, c(I_{p}))}{2^{k}} \geq \frac{{\mbox{dist}}(J, I_{p})}{2^{k}} + \frac{1}{2}.\]
For all $p \in \bf{T}$ with $|I_{p}| = 2^{kn}$, the $I_{p}$ are pairwise disjoint and contained in
$I_{\bf{T}}$. Therefore ${\mbox{dist}}(J, I_{p}) \geq {\mbox{dist}}(J, I_{{\bf{T}}})$ and 
$|I_{{\bf{T}}}|^{- \frac{1}{n}} \leq 2^{-k}$, which gives
\begin{equation} 
\left( 1 + \frac{{\mbox{dist}}\left(J, I_{p}\right)}{2^{k}} \right)^{-5n} \leq 
\left( 1 + \frac{{\mbox{dist}}\left(J, I_{{\bf{T}}}\right)}{\left|I_{{\bf{T}}} 
\right|^{\frac{1}{n}}} \right)^{-5n}. \label{k1est}\end{equation}
We will treat the remaining powers with the following lemma. 
\begin{lemma}
For $J \in \mathcal J$ such that $2^{kn} \leq 2^{n} |J|$,
\[ \sum_{\begin{subarray}{c} p \in \bf{T}\\ |I_{p}|= 2^{kn}\end{subarray}}\left( 1 + \frac{{\mathrm{dist}}
\left(J, I_{p}\right)}{2^{k}} \right)^{-5n} \leq C(n),\]
where $C(n)$ is independent of $J, k$ and ${\bf{T}}$. 
\end{lemma}

\begin{proof} We first observe that ${\mbox{dist}}(J, I_{p})$ and  ${\mbox{dist}}(c(J), I_{p})$ 
are of comparable size. The inequality ${\mbox{dist}}(J, I_{p}) \leq {\mbox{dist}}(c(J), I_{p})$ is clear.
 To see the other inequality, note that $|I_{p}| \leq 2^{n} |J| = |2J|$ implies that $I_{p}$ is 
 disjoint from $3J$, since $3J$ does not contain any $I_{p}$. Thus we have 
\begin{align*}
{\mbox{dist}}(I_{p}, c(J)) & \leq {\mbox{dist}}(I_{p}, J) + {\mbox{dist}}(\partial J, c(J)) \\
& \leq {\mbox{dist}}(I_{p}, J) + \frac{\sqrt{n}}{2} {\mbox{dist}}(I_{p}, J) \\
& \left(1 + \frac{\sqrt{n}}{2} \right) {\mbox{dist}}(I_{p}, J).
\end{align*}
Hence it suffices to replace ${\mbox{dist}}(J, I_{p})$ by ${\mbox{dist}}(c(J), I_{p})$. 
Let $\x_{0} = c(J)$ and decompose $\mathbb R^{n}$ as follows :
\[ \mathbb R^{n} = \bigcup_{m = 1}^{ \infty} \mathcal O_{m},\]
where
\begin{align*} \mathcal O_{1} &:= B(\x_{0}, 3 \sqrt{n} 2^{k}), \\
\mathcal O_{m} &:= B(\x_{0}, 3m \sqrt{n} 2^{k}) \backslash B(\x_{0}, 3(m-1)\sqrt{n} 2^{k}).
\end{align*}
Let 
\begin{align*}
S_{1} & := \left\{p \in {\bf{T}} \,:\, |I_{p}| = 2^{kn}, I_{p} \cap B \left(\x_{0}, 
3 \sqrt{n} 2^{k} \right) \ne \emptyset \right\}, \\
S_{m} & := \left\{p \in {\bf{T}} \,:\, |I_{p}| = 2^{kn}, I_{p} \cap \mathcal O_{m} \ne \emptyset, I_{p} 
\cap \left( \cup_{i=1}^{m-1} \mathcal O_{i}\right) = \emptyset \right\}, \quad m \geq 2. 
\end{align*}
Since the diameter of $I_{p}$ is $\sqrt{n}2^{k}$, $I_{p}$ will not intersect three annuli, 
so each $p$ in the sum is contained in exactly one $S_{m}$.  In order to estimate the number of 
tiles in $S_{m}$, we consider the volume of the corresponding annulus $\mathcal O_{m}$. Now, 
\[ {\mbox{volume}}(\mathcal O_{m}) = (3 \sqrt{n})^{n} 2^{kn} \left(m^{n} - (m-1)^{n} \right) 
= C_{n} 2^{kn} m^{n-1}.\]
Since the $I_{p}$-s are disjoint, there are $C_{n}m^{n-1}$ cubes $I_{p}$ of size $2^{kn}$ in the set 
$S_{m}$. Also for $p \in S_{m}$, 
\[ 3(m-1)\sqrt{n} 2^{k} \leq {\mbox{dist}}(\x_{0},I_{p}) \leq 3m \sqrt{n} 2^{k}.   \]
Thus, 
\[\sum_{ \begin{subarray}{c}p \in \bf{T}\\ |I_{p}| = 2^{kn} \end{subarray}}
\frac{1}{\left( 1 + \frac{{\mbox{dist}}\left( \x_{0}, I_{p}\right)}{2^{k}} \right)^{5n}} 
\leq \sum_{m=1}^{\infty} \frac{m^{n-1}}{(1+m)^{5n}} \leq \sum_{m=1}^{\infty} \frac{1}{m^{6}} < \infty.
 \]\end{proof}

Using (\ref{k1est}) and the lemma, we have that $\mathcal K_{1}$ is bounded by 
\begin{align*}
& C \epsilon \mu \sum_{J \in \mathcal J} \sum_{kn = - \infty }^{\log_{2} 2^{n}|J|}2^{kn}\left( 1 + 
\frac{{\mbox{dist}}\left(J, I_{{\bf{T}}}\right)}{\left|I_{{\bf{T}}} \right|^{\frac{1}{n}}} \right)^{-5n} 
\\ 
\leq & C \epsilon \mu \sum_{J \in \mathcal J} |J|\left( 1 + \frac{{\mbox{dist}}
\left(J, I_{{\bf{T}}}\right)}{\left|I_{{\bf{T}}} \right|^{\frac{1}{n}}} \right)^{-5n}\\
\leq & C \epsilon \mu \sum_{J \in \mathcal J} \int_{J} 
\left( 1 + \frac{\left| {\bf{x}} - c(I_{{\bf{T}}})\right|}{\left|I_{{\bf{T}}} \right|^{\frac{1}{n}}} \right)^{-5n}
\, d{\bf{x}}\\
\leq & C \epsilon \mu |I_{{\bf{T}}}|.
\end{align*}
This completes the estimate of $\mathcal K_{1}$. 

Now we consider $\mathcal K_{2}$ defined by (\ref{k2}). We can assume that
the summation runs only over those $J \in \mathcal J$ for which there
exists a $p \in {\bf{T}}$ with $2^n |J| < |I_{p}|$. Then we have $J \subset 3
I_{{\bf{T}}}$ and $2^n |J| < |I_{{\bf{T}}}|$ for all $J$ occurring in the sum. 

Let us fix a dyadic cube $J \in \mathcal J$ and observe that the set
\begin{displaymath}
G_{J} =J \cap \bigcup_{p \in {\bf{T}}\,:\, |I_{p}| > 2^n |J|} E_{2p}
\end{displaymath}
has measure at most $C \mu |J|$. To see this, let $J'$ be the unique dyadic cube
which contains $J$ and $|J'| = 2^{n} |J| < |I_{{\bf{T}}}|$. By maximality
of $J$, $3J'$ contains $I_{p_{0}}$ for some $p_{0} \in {\bf{T}}$. 
There are two cases to consider.  Case (a):  $I_{p_0}$ is the dyadic cube that is formed from taking 
the unique double of 
each side of $J'$ which is also dyadic.  In this case $|I_{p_0}|=2^n|J'|$ and we
set $p_0=p' < I_{{\bf{T}}} \times \omega_{{\bf{T}}}$.  Case (b): $I_{p_0}$ is contained in one of the dyadic cubes of size 
$|J'|$ contained in $3J'$.  Since  
$|J'|=2^n|J| <|\itt|$, the dyadic cube which contains $I_{p_0}$ is contained in $\itt$.
In this case there exists a tile $p'$ with $|I_{p'}|=|J'|$ so that
$I_{p_0} \subset I_{p'} \subset \itt$.  In both cases we have a tile $p'$ such that 
$p_o<p'< I_{{\bf{T}}} \times \omega_{{\bf{T}}}$ and $|\Om_{p'}|$ is either $2^{-n}|J|^{-1}$ or 
$2^{-2n}|J|^{-1}$.
We claim that 
\begin{displaymath}
\bigcup_{p \in {\bf{T}}\,:\, |I_{p}| > 2^n |J|} E_{2p} \subset E \cap \N^{-1}[\om_{p'}].
\end{displaymath}
To see this, let us choose $p \in {\bf{T}}$ such that $|I_{p}| > 2^n |J|$. Then
$|\omega_{p}| < 2^{-n}|J|^{-1}$, which means 
$|\omega_{p}| < |\omega_{p'}|$.
But $\omega_{{\bf{T}}} \subset \omega_{p} \cap \omega_{p'}$, which leads us to
conclude
$\omega_{p} \subset \omega_{p'}$.
Recalling that 
\begin{displaymath}
E_{2p} = \left\{{\bf{x}}\,:\, {\bf{N}}({\bf{x}}) \in \omega_{p(2^n)} \right\} \cap E 
\end{displaymath}
now completes the proof of the claim. \\
The above claim implies that $G_{J} \subset J' \cap E \cap \N^{-1}[\om_{p'}]$. Therefore, 
\begin{displaymath}
|G_{J}| \leq |J' \cap E \cap \N^{-1}\left[\om_{p'}\right]| = \intl_{E \cap \N^{-1}\left[ \om_{p'}
\right]} 1_{J'}({\bf{x}})\,
 d{\bf{x}}.
\end{displaymath}
Since 
\begin{displaymath}
1_{J'}({\bf{x}}) \leq C  \left(1 + \frac{\left|{\bf{x}} - c(I_{p'})\right|}{\left|I_{p'}\right|^
{\frac{1}{n}}}\right)^{- \nu},
\end{displaymath}
and ${\mbox{mass}}(\left\{p\right\}) \leq \mu$, we get $|G_{J}| \leq C
\mu |J|$. 

Let ${\bf{T}}_{2}$ be the $2^{n}$-tree of all $p \in {\bf{T}}$ such that $\omega_{{\bf{T}}(2^{n})}
\subset \omega_{p(2^{n})}$ and let ${\bf{T}}_{1} = {\bf{T}} \backslash {\bf{T}}_{2}$. Define for $j
= 1, 2,$
\begin{displaymath}
F_{jJ} := \sum_{p \in {\bf{T}}_{j}\,:\,|I_{P}| > 2^n |J|} \alpha_{P} \langle f,
\phi_{p} \rangle \psi_{p}^{{\bf{N}}} 1_{E_{2p}}.
\end{displaymath}
First we consider $F_{1J}$. We have 
\begin{align*}
|F_{1J}({\bf{x}})| & \leq \sum_{p \in {\bf{T}}_{1} \,:\, |I_{p}| > 2^n |J|} |\langle
f, \phi_{p} \rangle| |\psi_{p}^{{\bf{N}}} ({\bf{x}})| 1_{E_{2p}} ({\bf{x}})\\
& \leq C \epsilon \sum_{p \in {\bf{T}}_{1} \,:\, |I_{p}| > 2^n|J|} \left(1 +
\frac{|{\bf{x}} - c(I_{p})|}{|I_{p}|^{\frac{1}{n}}} \right)^{- \nu} 1_{E_{2p}}({\bf{x}}).
\end{align*}
We will sum the expression on the right hand side of the above inequality in two steps. First let
us construct
\begin{align*} \mathcal I &:= \left\{ \om \,:\, \textup{there exists $p \in \bf{T}_{1}$ such that $|I_{p}| >
2^n |J|$ and $\om_{p(2^{n})} = \om$} \right\}, \textup{ and }\\
\mathcal P_{\om} &:= \left\{ p \in {\bf{T}}_{1} \,:\, |I_{p}|>2^n |J|, \om_{p(2^n)} = \om \right\},
 \textup{ for } \om \in \mathcal I.
\end{align*}
This means that
the sum estimating $F_{1J}$ may be written as 
\[ \sum_{\om \in \mathcal I} \sum_{p \in \mathcal P_{\om}} \left(1 +
\frac{|{\bf{x}} - c(I_{p})|}{|I_{p}|^{\frac{1}{n}}} \right)^{- \nu} 1_{E_{\om}}({\bf{x}}),\]
where \[ E_{\om} := E \cap \left\{ \x \,:\, \N(\x) \in \om \right\}.\]
Now note that for $p \in {\bf{T}}_{1}$, the semitiles $I_{p} \times \om_{p(2^n)}$ are disjoint. In
particular, for $p, p' \in \mathcal P_{\om}$, $p \ne p'$, one has $I_{p} \cap I_{p'} = \emptyset$.
Therefore, 
\begin{displaymath}
\sum_{ p \in \mathcal P_{\om}}  \left(1 +
\frac{|{\bf{x}} - c(I_{p})|}{|I_{p}|^{\frac{1}{n}}} \right)^{- \nu} \leq
C.
\end{displaymath} 
The proof of this fact is similar to that of Lemma 5.
This implies that  
\begin{displaymath}
|F_{1J}({\bf{x}})| \leq C \epsilon \sum_{\om \in \mathcal I} 1_{E_{\om}}({\bf{x}}) = C \epsilon 
1_{\bigcup E_{\om}}({\bf{x}}).
\end{displaymath}
Here we have used the fact that the $\om$-s in $\mathcal I$ are
disjoint. This  yields
\begin{displaymath}
\left|\left|F_{1J}({\bf{x}}) \right|\right|_{L^{1}\left(J\right)} \leq C
\epsilon \int_{J}1_{\bigcup E_{\om}}({\bf{x}}) \, d{\bf{x}} = C \epsilon
 |G_{J}| \leq C
 \epsilon \mu |J|,
\end{displaymath}
This estimate, summed over the disjoint $J \subset 3I_{{\bf{T}}}$ yields the
desired bound. 

To complete the proof of (\ref{45}) we estimate
$F_{2J}({\bf{x}})$. Fix ${\bf{x}}$ and assume that $F_{2J}({\bf{x}})$ is not
zero. Since the cubes $\omega_{p(2^n)}$ with $p \in {\bf{T}}_{2}$ are all nested
and 
$E_{2p} = \left\{ {\bf{x}} \,:\, {\bf{N}}({\bf{x}}) \in \omega_{p(2^n)}
\right\} \cap E$, 
there is a largest cube $\omega_{+}$ of the form
$\omega_{p}$ with $p \in {\bf{T}}_{2},\; {\bf{x}} \in E_{2p}$ and $|I_{p}| >
2^n |J|$. Similarly there is a smallest cube which we call $\om_{s}$ satisfying the above
properties. Let us define $\om_{-} = \om_{s(2^n)}$. Then $ {\bf{x}} \in E_{2p}$ for some 
$p \in {\bf{T}}$ with $|I_{p}| > 2^n|J|$ if
and only if  $|\omega_{-}| < | \omega_{p} | \leq | \omega_{+} |$. Fix
$ {\pmb{\xi}}_{0} \in \omega_{{\bf{T}}}$. We can now write $F_{2J}({\bf{x}})$
as 
\begin{displaymath}
F_{2J}({\bf{x}}) = \sum_{\begin{subarray}{c}p \in {\bf{T}}_{2}\\ |\omega_{-}| < |\omega_{p}|
\leq |\omega_{+}| \end{subarray}} \alpha_{p} \langle f, \phi_{p} \rangle
\psi_{p}^{{\bf{N}}} ({\bf{x}}),
\end{displaymath}
which may be decomposed as  
\begin{multline*}
\sum_{\begin{subarray}{c}p \in {\bf{T}}_{2}\\ |\omega_{-}| < |\omega_{p}|
\leq |\omega_{+}| \end{subarray}} \alpha_{p}  \langle f, \phi_{p}
\rangle \left(e^{ 2 \pi i {\pmb{\xi}}_{0} \cdot \left( \cdot \right)} K (
\cdot) \ast \phi_{p} (\cdot) \right) ({\bf{x}})   \\
+\left( \left[\left(e^{2 \pi i {\bf{N}} \left({\bf{x}}\right) \cdot \left( \cdot \right)} - e^{2 \pi i
{\pmb{\xi}}_{0} \cdot \left( \cdot \right)} \right) K( \cdot) \right]  
\ast \sum_{\begin{subarray}{c}p \in {\bf{T}}_{2}\\ |\omega_{-}| < |\omega_{p}|
\leq |\omega_{+}| \end{subarray}} \alpha_{p} \langle f, \phi_{p}
\rangle \phi_{p} (\cdot ) \right) ({\bf{x}})
\end{multline*}
\begin{multline*}
= \sum_{p \in {\bf{T}}_{2}} \alpha_{p} \langle f, \phi_{p} \rangle \left(
\psi_{p}^{{\pmb{\xi}}_{0}} \ast 
\left( M_{c\left(\omega_{+} \right)} D_{\frac{1}{6} |\omega_{+}|^{- \frac{1}{n}}}^{1}
\phi - M_{c\left(\omega_{-} \right)} D_{\frac{1}{6} |\omega_{-}|^{- \frac{1}{n}}}^{1}
\phi \right) \right) ({\bf{x}})\\
+ \Bigl\{ \left[\left(e^{2 \pi i {\bf{N}} \left({\bf{x}}\right) \cdot \left( \cdot \right)} - e^{2 \pi i
{\pmb{\xi}}_{0} \cdot \left( \cdot \right)} \right) K( \cdot) \right]
 \\
\ast \Bigl[ \sum_{p \in {\bf{T}}_{2}} \alpha_{p} \langle f, \phi_{p} \rangle \left(
\phi_{p} \ast 
\left( M_{c\left(\omega_{+} \right)} D_{\frac{1}{6} |\omega_{+}|^{- \frac{1}{n}}}^{1}
\phi - M_{c\left(\omega_{-} \right)} D_{\frac{1}{6} |\omega_{-}|^{- \frac{1}{n}}}^{1}
\phi \right) \right) \Bigr]\Bigr\}  ({\bf{x}}).
\end{multline*} 
We claim that the last equality follows from the geometry of the supports of the Fourier transforms of the two
convolving functions. More specifically, $\widehat \phi_p$ is supported on $\frac{1}{5} \om_{p(1)}$
while 
$$
\left( M_{c\left(\omega_{\pm} \right)} D_{\frac{1}{6} |\omega_{\pm}|^{- \frac{1}{n}}}^{1}\phi \right)
\widehat{} \, \, ({\pmb {\xi}})= 
\begin{cases}
1  & \textup{if } {\pmb {\xi}} \in \om_{\pm} \\
0  &  \textup{if } {\pmb {\xi}} \notin \left(1 + \frac{1}{5}\right)\om_{\pm}.
\end{cases} 
$$
Therefore
\begin{multline*}
\left(M_{c\left(\omega_{+} \right)} D_{\frac{1}{6} |\omega_{+}|^{- \frac{1}{n}}}^{1}
\phi - M_{c\left(\omega_{-} \right)} D_{\frac{1}{6} |\omega_{-}|^{- \frac{1}{n}}}^{1}
\phi\right)\widehat{} \, \, ({\pmb{\xi}})= \\
\begin{cases}
1 &  \textup{ if } {\pmb {\xi}} \in \om_{+} \setminus 
\left(1 + \frac{1}{5}\right)\om_{-}\\
0 & \textup{if } {\pmb {\xi}} \in \om_{-} \cup \left( \left (1 + \frac{1}{5}\right )\om_{+}\right)^c.
\end{cases}
\end{multline*}
For those $p \in \mathbf T_2$ such that $|\om_-| < |\om_p| \le |\om_+|$, we have
 $$\frac{1}{5} \om_{p(1)} \subset \om_{+} \setminus 
\left(1 + \frac{1}{5}\right)\om_{-}.$$  Conversely, for $p \in \mathbf T_2$ such that
$|\om_p| > |\om_{+}|$ or $|\om_p| \le |\om_{-}|$, we have 
$$\frac{1}{5} \om_{p(1)} \subset
\om_{-} \cup \left( \left (1 + \frac{1}{5}\right )\om_{+}\right)^c.$$
 This concludes the proof of the claim.
 
 The expression for $F_{2J}$ may therefore be written as 
\begin{multline}
\left(G_{1} \ast
\left( M_{c\left(\omega_{+} \right)} D_{\frac{1}{6} |\omega_{+}|^{- \frac{1}{n}}}^{1}
\phi - M_{c\left(\omega_{-} \right)} D_{\frac{1}{6}\omega_{-}|^{- \frac{1}{n}}}^{1}
\phi \right) \right)({\bf{x}})  \\
+ \biggl[ \! \left(e^{2 \pi i {\bf{N}} \left({\bf{x}}\right) \cdot \left( \cdot \right)} - e^{2 \pi i
{\pmb{\xi}}_{0} \cdot \left( \cdot \right)} \right) K( \cdot) \ast
G_{2}(\cdot) \\
\ast \left( M_{c\left(\omega_{+} \right)} D_{\frac{1}{6} |\omega_{+}|^{- \frac{1}{n}}}^{1}
\phi - M_{c\left(\omega_{-} \right)} D_{\frac{1}{6} |\omega_{-}|^{- \frac{1}{n}}}^{1}
\phi \right) \! \biggr]  ({\bf{x}}), \label{sum}
\end{multline}
where 
\begin{align*}
G_{1}({\bf{x}}) &= \sum_{p \in {\bf{T}}_{2}} \alpha_{p} \langle f, \phi_{p}
\rangle \psi_{p}^{{\pmb{\xi}}_{0}} ({\bf{x}}) \\
& = \sum_{p \in {\bf{T}}_{2}} \alpha_{p} \langle f, \phi_{p}\rangle \left(
e^{2 \pi i {\pmb{\xi}}_{0} \cdot \left( \cdot \right)} K(\cdot) \ast \phi_{p}
(\cdot) \right)({\bf{x}}),\\
G_{2}({\bf{x}}) &= \sum_{p \in {\bf{T}}_{2}} \alpha_{p} \langle f, \phi_{p}
\rangle \phi_{p}.
\end{align*}

\begin{claim}
\begin{equation}\label{claim}
|F_{2J}({\bf{x}})| \leq C \left( \sup_{J \subset I} \frac{1}{|I|}
 \int_{I} |G_{1}({\bf{z}})|\,d{\bf{z}} +\sup_{J \subset I} \frac{1}{|I|}
 \int_{I} |G_{2}({\bf{z}})|\,d{\bf{z}} \right),  
\end{equation}
where the suprema above are taken over all cubes $I$ containing $J$.
\end{claim}

The proof of the claim is given in the next section. One should
recognize the claim as a slight variant of the classical inequality 
\begin{displaymath}
T^{*}\bar{f} \lesssim M(T\bar{f}) + M({\bar{f}}),
\end{displaymath}
where 
\begin{displaymath}
Tg = \left(e^{2 \pi i {\pmb{\xi}}_{0} \cdot \left( \cdot \right)} K(\cdot)
\right) \ast g, \quad \bar{f} = G_{2}(x), 
\end{displaymath}
$T^*$ is the maximal operator corresponding to $T$, 
and $M$ denotes the Hardy-Littlewood maximal function. In the rest of this
section we show how the proof of Theorem 1 may be
completed using (\ref{claim}). 

We observe that the right hand side of the above expression is constant
on $J$ and that $F_{2J}1_{J}$ is supported on $G_{J}$, which is of
measure less than or equal to $C \mu |J|$. Hence, 
\begin{align*}
& \sum_{J \in \mathcal J} ||F_{2J}||_{L^{1}\left(J\right)}  \\
 \leq &   C \mu
\sum_{J \in \mathcal J \,:\, J \subset 3 I_{{\pmb{T}}}} |J| \left( \sup_{J \subset I} \frac{1}{|I|}
 \int_{I} |G_{1}({\bf{z}})|\,d{\bf{z}} +\sup_{J \subset I} \frac{1}{|I|}
 \int_{I} |G_{2}({\bf{z}})|\,d{\bf{z}} \right)  \\
 \leq  &C \mu \left( \biggl\| M \biggl( \sum_{p \in {\bf{T}}_{2}} \alpha_{p}
 \langle f, \phi_{p} \rangle \psi_{p}^{{\pmb{\xi}}_{0}} \biggr)
 \biggr\|_{L^{1}\left(3 I_{{\bf{T}}}\right)} \! \! \! \! +\biggl\| M \biggl( \sum_{p \in {\bf{T}}_{2}} 
 \alpha_{p}
 \langle f, \phi_{p} \rangle \phi_{p} \biggr)
 \biggr\|_{L^{1}\left(3 I_{{\bf{T}}}\right)} \right)\\
 \leq & C \mu |I_{{\bf{T}}}|^{\frac{1}{2}} \left( \biggl\| \sum_{p \in {\bf{T}}_{2}} \alpha_{p}
 \langle f, \phi_{p} \rangle \psi_{p}^{{\pmb{\xi}}_{0}}
 \biggr\|_{L^{2}\left( \mathbb R^{n}\right)} \! \! \! \! + \biggl\| \sum_{p \in {\bf{T}}_{2}} \alpha_{p}
 \langle f, \phi_{p} \rangle \phi_{p}
 \biggr\|_{L^{2}\left(\mathbb R^{n} \right)} \right).
\end{align*}
Here we have used the $L^{2}$ boundedness of the
Hardy-Littlewood maximal function $M$. We would now like to show that $||G_{1}||_{L^{2}}$ and
$||G_{2}||_{L^{2}}$ are bounded above by a constant multiple of $\epsilon |I_{T}|^{\frac{1}{2}}$.    
\begin{align*}
||G_{1}||_{L^{2}}^{2} & = \sum_{p, p' \in {\bf{T}}_{2}} \alpha_{p}
  \alpha_{p'} \langle f, \phi_{p} \rangle \; \overline{\langle f, \phi_{p'}
  \rangle} \; \langle \psi_{p}^{{\pmb{\xi}}_{0}},
  \psi_{p'}^{{\pmb{\xi}}_{0}} \rangle \\
& = \sum_{\begin{subarray}{c}p, p' \in {\bf{T}}_{2} \\ \om_p \ne \om_p' \end{subarray}} \alpha_{p}
  \alpha_{p'} \langle f, \phi_{p} \rangle \; \overline{\langle f, \phi_{p'}
  \rangle} \; \langle m( \cdot - {\pmb{\xi}}_{0}) \widehat{\phi_{p}}, m(
  \cdot - {\pmb{\xi}}_{0}) \widehat{\phi_{p'}} \rangle \\
& \qquad  \qquad \qquad \qquad  + \sum_{\begin{subarray}{c} p, p' \in {\bf{T}}_{2} \\ \om_p = \om_p' 
\end{subarray}}\alpha_{p}
  \alpha_{p'} \langle f, \phi_{p} \rangle \; \overline{\langle f, \phi_{p'}
  \rangle} \; \langle \psi_{p}^{{\pmb{\xi}}_{0}},
  \psi_{p'}^{{\pmb{\xi}}_{0}} \rangle.
\end{align*}
Similarly, 
\begin{align*}
||G_{2}||_{L^{2}}^{2} & = \sum_{p, p' \in {\bf{T}}_{2}} \alpha_{p}
  \alpha_{p'} \langle f, \phi_{p} \rangle \; \overline{\langle f, \phi_{p'}
  \rangle} \; \langle {\phi}_{p}, {\phi}_{p'} \rangle \\
& = \sum_{ \begin{subarray}{c}p, p' \in {\bf{T}}_{2} \\ \om_p \ne \om_p' \end{subarray}} \alpha_{p}
  \alpha_{p'} \langle f, \phi_{p} \rangle \; \overline{\langle f, \phi_{p'}
  \rangle} \; \langle \widehat{\phi_{p}}, \widehat{\phi_{p'}} \rangle \\
& \qquad \qquad \qquad \qquad + \sum_{ \begin{subarray}{c}p, p' \in
  {\bf{T}}_{2} \\ \om_p = \om_p' \end{subarray}} \alpha_{p}
  \alpha_{p'} \langle f, \phi_{p} \rangle \; \overline{\langle f, \phi_{p'}
  \rangle} \; \langle {\phi}_{p}, {\phi}_{p'} \rangle.
\end{align*}
But $\langle \widehat{\phi_{p}}, \widehat{\phi_{p'}}\rangle = 0$ 
for $p, p' \in {\bf{T}}_{2}, \om_p \ne \om_{p'}$, since $\widehat{\phi_{p}}$ and
$\widehat{\phi_{p'}}$ have disjoint supports in this case. Therefore we only need to consider
the second sum in the expressions for $||G_{1}||_{L^{2}}^{2}$ and
$||G_{2}||_{L^{2}}^{2}$. Our pointwise bounds imply that 
\begin{multline*}
 \langle {\phi}_{p}, {\phi}_{p'} \rangle,\langle \psi_{p}^{{\pmb{\xi}}_{0}},
  \psi_{p'}^{{\pmb{\xi}}_{0}} \rangle \lesssim \\
  |I_{p}|^{- \frac{1}{2}}
  |I_{p'}|^{- \frac{1}{2}} \int \left( 1 + \frac{ \left| {\bf{x}} - c(I_{p})\right|}
  {|I_{p}|^{\frac{1}{n}}}\right)^{-\nu} \left( 1 + \frac{ \left| {\bf{x}} - c(I_{p'})\right|}
  {|I_{p'}|^{\frac{1}{n}}}\right)^{-\nu} \, d{\bf{x}},
\end{multline*}
so it is enough to estimate the right hand side above for $p, p'$ satisfying
$\om_p = \om_p'$.  Upon simplification this reduces to
\begin{align*}
& |I_{p}|^{-1}
\int_{\mathbb R^{n}} \left(1 + \frac{|{\bf{x}} -
c(I_{p})|}{|I_{p}|^{\frac{1}{n}}} \right)^{- \nu} \left(1 + \frac{|{\bf{x}} - c(I_{p'})|}
{|I_{p'}|^{\frac{1}{n}}}\right)^{- \nu}\, d{\bf{x}}\\      
= &\int_{\mathbb R^{n}} \left(1 + |{\bf{y}}| \right)^{- \nu} \left(1 +
\left|{\bf{y}} + \frac{c(I_{p}) - c(I_{p})}{|I_{p}|^{\frac{1}{n}}} \right| \right)^{-\nu}
\, d{\bf{y}}\\
 \lesssim & \left(1 + \left| \frac{c(I_{p}) - c(I_{p'})}{|I_{p}|^{\frac{1}{n}}}
\right|\right)^{- \nu}.
\end{align*}
With this estimate, the proof that $||G_{1}||_{L^{2}}$ and $||G_{2}||_{L^{2}}$ are less than 
$C \epsilon |I_{\bf{T}}|^{\frac{1}{2}}$ is similar to an argument outlined in the proof of Lemma 2. 
One needs to follow the proof of the estimate of the first term of (\ref{65}) to complete 
the proof of Lemma 3, given the claim.


\section{Proof of Claim 1}
\setcounter{equation}0

Let us estimate the first term in the expression (\ref{sum}). We denote
by $\tilde c$ any one of the two constants $\frac{1}{6} |\omega_{+}|^{-
\frac{1}{n}}$ or $\frac{1}{6} |\omega_{-}|^{- \frac{1}{n}}$. By translation invariance,
let $I_{-}$ and
$I_{+}$ be the unique dyadic cubes
of the form 
$$
I_{\pm}= \prod_{j=1}^n \, [0, |\om_{\pm}|^{-\frac{1}{n}}).
$$
For a dyadic cube $I=\prod_{j=1}^n \, [0,2^k)$  and 
$\mathbf r = (r_1,r_2, \cdots,r_n) \in \mathbb Z^n$, $I + \mathbf r$ will denote the unique dyadic cube
of the form
$$
[r_1 2^k, (r_1 +1) 2^k) \times 
\cdots \times [r_n 2^k, (r_n +1) 2^k).
$$
Now,
\begin{align*}
& \int |G_{1}({\bf{y}})| {\tilde{c}}^{-n} \left| \phi \left(\frac{ {\bf{x}}
- {\bf{y}}}{\tilde{c}} \right) \right|\, d{\bf{y}} \\
\leq & {\tilde{c}}^{-n}
\intl_{|{\bf{x}} - {\bf{y}}| \leq \tilde{c}} |G_{1}({\bf{y}})| \,
d{\bf{y}} + {\tilde{c}}^{-n} \sum_{j \geq 1} 2^{-j \nu} \! \! \! \! \intl_{|{\bf{x}}
- {\bf{y}}| \sim \tilde{c} 2^{j}} |G_{1}({\bf{y}})|\, d{\bf{y}} \\
\leq & {\tilde{c}}^{-n}
\intl_{ \begin{subarray}{c}|x_{i} - y_{i}| \leq \tilde{c} \\
1\leq i \leq n \end{subarray}} |G_{1}({\bf{y}})| \,
d{\bf{y}}  
+ \sum_{j \geq 1} 2^{-j} \frac{1}{\left(\tilde{c} 2^{j}\right)^{n}}
\intl_{\begin{subarray}{c}|x_{i} - y_{i}| \leq \tilde{c}2^{j}\\ 1 \leq i \leq n\end{subarray}}
|G_{1}({\bf{y}})|\, d{\bf{y}}\\
 \lesssim & \sup_{I\,:\,J \subset I} \frac{1}{|I|} \int_{I}
|G_{1}({\bf{z}})|\, d{\bf{z}}.
\end{align*}
Here we have used the fact that since ${\bf{x}} \in J$ and $ |J| <
|I_{+}| < |I_{-}|$, we have $J \subset \left\{ {\bf{y}}\,:\, |x_{i} -
y_{i} | \leq \tilde{c} {\mbox{ for all }} i, 1 \leq i \leq n \right\}$. \\

We denote by $\mathcal B$ the second term in (\ref{sum}).
{\allowdisplaybreaks 
\begin{multline*}
\mathcal B = \int_{({\bf{y,z}}) \in \mathbb R^{n} \times \mathbb
R^{n}} \left(e^{2 \pi i {\bf{N}} \left( {\bf{x}} \right) \cdot {\bf{y}}} -
e^{2 \pi i {\pmb{\xi}}_{0} \cdot {\bf{y}} } \right) K({\bf{y}})  \\ \times \Bigl[e^{2 \pi i c\left(
\omega_{+} \right) \cdot \left( {\bf{z}} - {\bf{y}} \right)} \left( \frac{1}{6}
|\omega_{+}|^{- \frac{1}{n}} \right)^{-n} \phi \left(\frac{ {\bf{z}} -
{\bf{y}}}{\frac{1}{6} |\omega_{+}|^{-\frac{1}{n}} } \right) \\
- e^{2 \pi i c\left(
\omega_{-} \right) \cdot \left( {\bf{z}} - {\bf{y}} \right)} \left( \frac{1}{6}
|\omega_{-}|^{- \frac{1}{n}} \right)^{-n} \phi \left(\frac{ {\bf{z}} -
{\bf{y}}}{\frac{1}{6} |\omega_{-}|^{-\frac{1}{n}} } \right)\Bigr]
G_{2}({\bf{x}} - {\bf{z}}) \, d{\bf{y}} d{\bf{z}}.
\end{multline*}}
To estimate
$\mathcal B$ we write it as 
\begin{displaymath}
\mathcal B = \mathcal B_{1} - \mathcal B_{2} + \mathcal B_{3} + \mathcal B_{4},
\end{displaymath}
where 
\begin{multline*}
\mathcal B_{1} := \int_{{\bf{z}} \in I_{+}} \int_{{\bf{y}} \in \mathbb R^{n}}
\left(e^{2 \pi i {\bf{N}} \left( {\bf{x}} \right) \cdot {\bf{y}}} -
e^{2 \pi i {\pmb{\xi}}_{0} \cdot {\bf{y}} } \right) K({\bf{y}})
e^{2 \pi i c\left(
\omega_{+} \right) \cdot \left( {\bf{z}} - {\bf{y}} \right)} 
\\ \times \left( \frac{1}{6}
|\omega_{+}|^{- \frac{1}{n}} \right)^{-n}  \phi \left(\frac{ {\bf{z}} -
{\bf{y}}}{\frac{1}{6} |\omega_{+}|^{-\frac{1}{n}} } \right) G_{2}({\bf{x}}
-{{\bf{z}}}) \, d{\bf{y}} \, d{\bf{z}}, 
\end{multline*}
\begin{multline*}
\mathcal B_{2} := \int_{{\bf{z}} \in I_{-}} \int_{{\bf{y}} \in \mathbb R^{n}}
\left(e^{2 \pi i {\bf{N}} \left( {\bf{x}} \right) \cdot {\bf{y}}} -
e^{2 \pi i {\pmb{\xi}}_{0} \cdot {\bf{y}} } \right) K({\bf{y}})
e^{2 \pi i c\left(
\omega_{-} \right) \cdot \left( {\bf{z}} - {\bf{y}} \right)} 
\\ \times \left( \frac{1}{6}
|\omega_{-}|^{- \frac{1}{n}} \right)^{-n}   \phi \left(\frac{ {\bf{z}} -
{\bf{y}}}{\frac{1}{6} |\omega_{-}|^{-\frac{1}{n}} } \right) G_{2}({\bf{x}}
-{{\bf{z}}}) \, d{\bf{y}} \, d{\bf{z}}, 
\end{multline*}
\begin{multline*}
\mathcal B_{3} := \sum_{{\bf{r}} \in \mathbb Z^{n} \setminus
\left\{{\bf{0}} \right\}} \int_{{\bf{z}} \in I_{-} + {\bf{r}}}
\int_{{\bf{y}} \in \mathbb R^{n}} \left(e^{2 \pi i {\bf{N}} \left( {\bf{x}} \right) \cdot {\bf{y}}} -
e^{2 \pi i {\pmb{\xi}}_{0} \cdot {\bf{y}} } \right) K({\bf{y}})  \\ \times \Bigl[e^{2 \pi i c\left(
\omega_{+} \right) \cdot \left( {\bf{z}} - {\bf{y}} \right)} \left( \frac{1}{6}
|\omega_{+}|^{- \frac{1}{n}} \right)^{-n} \phi \left(\frac{ {\bf{z}} -
{\bf{y}}}{\frac{1}{6} |\omega_{+}|^{-\frac{1}{n}} } \right) \\
- e^{2 \pi i c\left(
\omega_{-} \right) \cdot \left( {\bf{z}} - {\bf{y}} \right)} \left( \frac{1}{6}
|\omega_{-}|^{- \frac{1}{n}} \right)^{-n} \phi \left(\frac{ {\bf{z}} -
{\bf{y}}}{\frac{1}{6} |\omega_{-}|^{-\frac{1}{n}} } \right)\Bigr]
G_{2}({\bf{x}} - {\bf{z}}) \, d{\bf{y}} d{\bf{z}},
\end{multline*}
and 
\begin{multline*}
\mathcal B_{4} := \int_{{\bf{z}} \in I_{-} \setminus I_{+}}
\int_{{\bf{y}} \in \mathbb R^{n}}
\left(e^{2 \pi i {\bf{N}} \left( {\bf{x}} \right) \cdot {\bf{y}}} -
e^{2 \pi i {\pmb{\xi}}_{0} \cdot {\bf{y}} } \right) K({\bf{y}})
e^{2 \pi i c\left(
\omega_{+} \right) \cdot \left( {\bf{z}} - {\bf{y}} \right)}  \\ \times
\left( \frac{1}{6}
|\omega_{+}|^{- \frac{1}{n}} \right)^{-n} \phi \left(\frac{ {\bf{z}} -
{\bf{y}}}{\frac{1}{6} |\omega_{+}|^{-\frac{1}{n}} } \right) G_{2}({\bf{x}}
-{{\bf{z}}}) \, d{\bf{y}} \, d{\bf{z}}. 
\end{multline*}

We have the following bound for $\mathcal B_{1}$:
{\allowdisplaybreaks 
\begin{align*} 
|\mathcal B_{1}| & \leq \left[\int_{{\bf{z}} \in I_{+}} \int_{{\bf{y}} \in
 I_{+}}  + \sum_{{\bf{m}} \in \mathbb Z^{n} \setminus \left\{0\right\}}
 \int_{{\bf{z}} \in I_{+}} \int_{{\bf{y}} \in I_{+} + {\bf{m}}} \right]\\
&\left| \left(e^{2 \pi i {\bf{N}}\left({\bf{x}}\right)\cdot {\bf{y}}} - e^{2 \pi i
 {\pmb{\xi}}_{0} \cdot {\bf{y}}} \right) K({\bf{y}}) \right|
 |\omega_{+}| \left| \phi \left( \frac{{\bf{z}} - {\bf{y}}}{\frac{1}{6}
 \left|\omega_{+}\right|^{ - \frac{1}{n}}} \right) \right|\,d{\bf{y}}
 |G_{2}({\bf{x}} - {\bf{z}})| \, d{\bf{z}}\\
& \lesssim \int_{{\bf{z}} \in I_{+}} |\omega_{+}| \int_{{\bf{y}} \in I_{+}}
 |{\bf{N}}({\bf{x}}) - {\pmb{\xi}}_{0}| |{\bf{y}}|^{1-n} d{\bf{y}}
 |G_{2}({\bf{x}} - {\bf{z}})| \, d{\bf{z}} \quad  \\ 
&  \qquad \qquad +\sum_{{\bf{m}} \in \mathbb Z^{n} \setminus
 \left\{0\right\}} \int_{{\bf{z}} \in I_{+}}|\omega_{+}|\int_{{\bf{y}}
 \in I_{+} + {\bf{m}}} |{\bf{y}}|^{-n} |{\bf{m}}|^{-\nu} \,
 d{\bf{y}}|G_{2}({\bf{x}} - {\bf{z}})| \, d{\bf{z}}\\
& \lesssim |\omega_{+}| |\omega_{-}|^{\frac{1}{n}} |I_{+}|^{\frac{1}{n}}
 \int_{{\bf{z}} \in I_{+}} |G_{2}({\bf{x}} - {\bf{z}})| \, d{\bf{z}}\quad 
 \\ 
&  \qquad \qquad +\sum_{{\bf{m}} \in \mathbb Z^{n} \setminus
 \left\{0\right\}} |\omega_{+}| \left(|{\bf{m}}|\, |I_{+}|^{\frac{1}{n}}
 \right)^{-n} |{\bf{m}}|^{- \nu} |I_{+}| \int_{{\bf{z}} \in I_{+}}
 |G_{2}({\bf{x}} - {\bf{z}})|\,d{\bf{z}} \\
& \lesssim \left[
 \left(\frac{|\omega_{-}|}{|\omega_{+}|}\right)^{\frac{1}{n}} + \sum_{{\bf{m}} \in \mathbb Z^{n} 
\setminus
 \left\{0\right\}} |{\bf{m}}|^{-n - \nu} \right] \sup_{I\,:\, J \subset
 I} \frac{1}{|I|} \int_{I} |G_{2}({\bf{z}})| \, d{\bf{z}} \\
& \lesssim \sup_{I \,:\, J \subset I} \frac{1}{|I|} \int_{I} |G_{2}({\bf{z}})| \, d{\bf{z}},
\end{align*}}
for large $\nu$. The treatment for $\mathcal B_{2}$ is similar. \\

Next we estimate $\mathcal B_{3}$.  Note that it suffices to consider $|\mathbf r| \geq 2$.
The case $|\bf r|< 2$ follows with a similar argument as in the treatment of $\mathcal B_1$.
\begin{multline*}
\mathcal B_{3} = \sum_{{\bf{r}} \in \mathbb Z^{n} \setminus
\left\{0\right\}} \int_{{\bf{z}} \in I_{-} + {\bf{r}}} \int_{{\pmb{\xi}}
\in \mathbb R^{n}} \left( m({\pmb{\xi}} - {\bf{N}}({\bf{x}})) - m({\pmb{\xi}} -
{\pmb{\xi}}_{0})\right)e^{2 \pi i {\bf{z}} \cdot {\pmb{\xi}}}   \\ \times
\left[ \widehat{\phi} \left( \frac{1}{6} |\omega_{+}|^{ - \frac{1}{n}} ({\pmb{\xi}}
- c(\omega_{+})) \right) - \widehat{\phi} \left( \frac{1}{6} |\omega_{-}|^{ - \frac{1}{n}} ({\pmb{\xi}}
- c(\omega_{-})) \right)\right]\, d{\pmb{\xi}} |G_{2}({\bf{x}} -
{\bf{z}})| \, d{\bf{z}}.
\end{multline*}
Let us write
\begin{displaymath}
e^{2 \pi i {\bf{z}} \cdot {\pmb{\xi}}} = |{\bf{z}}|^{-2 N} (L_{\pmb{\xi}})^{N}
\left( e^{2 \pi i {\bf{z}} \cdot {\pmb{\xi}}}\right),
\end{displaymath}
where $N$ is a large positive integer and $L_{\pmb{\xi}}$ is a suitably chosen differential operator. 
Then we can see $\mathcal B_{3}$ as
\begin{multline*}
\sum_{{\bf{r}} \in \mathbb Z^{n} \setminus
\left\{0\right\}} \int_{{\bf{z}} \in I_{-} + {\bf{r}}}  |{\bf{z}}|^{- 2N}
\int_{{\pmb{\xi}}\in \mathbb R^{n}} e^{2 \pi i {\bf{z}} \cdot {\pmb{\xi}}}
L_{{\pmb{\xi}}}^{N} \biggl[\left( m({\pmb{\xi}} - {\bf{N}}({\bf{x}})) - m({\pmb{\xi}} -
{\pmb{\xi}}_{0})\right) \\ \times
\left( \widehat{\phi} \left( \frac{1}{6} |\omega_{+}|^{ - \frac{1}{n}} ({\pmb{\xi}}
- c(\omega_{+})) \right) - \widehat{\phi} \left( \frac{1}{6} |\omega_{-}|^{ - \frac{1}{n}} ({\pmb{\xi}}
- c(\omega_{-})) \right)\right) \biggr]\, d{\pmb{\xi}} |G_{2}({\bf{x}} -
{\bf{z}})| \, d{\bf{z}}.
\end{multline*}
For simplicity, let us consider only those terms where
$\left(L_{\pmb{\xi}}\right)^{N}$ is applied to either one of the terms  
\begin{multline*}
\left( m({\pmb{\xi}} - {\bf{N}}({\bf{x}})) - m({\pmb{\xi}} -
{\pmb{\xi}}_{0})\right) {\mbox{ or }}\\ \left( \widehat{\phi} \left( \frac{1}{6} |\omega_{+}|^{ - \frac{1}{n}} 
({\pmb{\xi}}
- c(\omega_{+})) \right) - \widehat{\phi} \left( \frac{1}{6} |\omega_{-}|^{ - \frac{1}{n}} ({\pmb{\xi}}
- c(\omega_{-})) \right)\right).
\end{multline*} The analysis for the other terms is
similar. First let us look at the inner integral 
\begin{multline} 
\int_{{\pmb{\xi}} \in \mathbb R^{n}} e^{2 \pi i {\bf{z}}
\cdot {\pmb{\xi}}}
L_{{\pmb{\xi}}}^{N} \left[\left( m({\pmb{\xi}} - {\bf{N}}({\bf{x}})) - m({\pmb{\xi}} -
{\pmb{\xi}}_{0})\right) \right]  \\ \times
 \left[\left( \widehat{\phi} \left( \frac{1}{6} |\omega_{+}|^{ - \frac{1}{n}} ({\pmb{\xi}}
- c(\omega_{+})) \right) - \widehat{\phi} \left( \frac{1}{6} |\omega_{-}|^{ - \frac{1}{n}} ({\pmb{\xi}}
- c(\omega_{-})) \right)\right) \right]\, d{\pmb{\xi}}. \label{innerint}
\end{multline}
We observe that 
\begin{displaymath}
\left| L_{{\pmb{\xi}}}^{N} \left[\left( m({\pmb{\xi}} - {\bf{N}}({\bf{x}})) - m({\pmb{\xi}} -
{\pmb{\xi}}_{0})\right) \right] \right| \lesssim \left( |{\pmb{\xi}} -
{\bf{N}}({\bf{x}})|^{-2N} + |{\pmb{\xi}} - {\pmb{\xi}}_{0}|^{-2N} \right),
\end{displaymath}
and that the integrand is supported on $(1 + \frac{1}{5}) \omega_{+} \setminus 
\omega_{-}$. Also, given $ \omega_{+}$ and $\omega_{-}$, there exists a
unique sequence of nested intervals 
\begin{displaymath}
\omega_{-} \subset \omega_{p_{1}} \subset \omega_{p_{2}} \subset \cdots
\subset \omega_{p_{M}}=\omega_{+}, \quad |\omega_{p_{i+1}}| = 2^{n} |
\omega_{p_{i}}|.
\end{displaymath}
It is not difficult to see that if ${\pmb{\xi}}\in (1 + \frac{1}{5}) \omega_{p_{i+1}}
\setminus   \omega_{p_{i}}$ and ${\bf{N}}({\bf{x}}), {\pmb{\xi}}_{0} \in
\omega_{-} \subset \omega_{p_{i}}$, then 
\begin{displaymath}
|{\bf{N}}({\bf{x}}) - {\pmb{\xi}}|, |{\pmb{\xi}} - {\pmb{\xi}}_{0}| \gtrsim |
 \omega_{p_{i}}|^{\frac{1}{n}}, \quad 1 \leq i \leq M-1. 
\end{displaymath}
This implies that the expression in (\ref{innerint}) is bounded by a
constant multiple of 
\begin{displaymath}
\sum_{i = 1}^{M-1} |\omega_{p_{i}}| | \omega_{p_{i}} |^{ -\frac{2N}{n}}
\lesssim |\omega_{-}|^{1 - \frac{2N}{n}}. 
\end{displaymath}
Next we consider the inner integral 
\begin{multline} 
\int_{{\pmb{\xi}} \in \mathbb R^{n}} e^{2 \pi i {\bf{z}}
\cdot {\pmb{\xi}}}
\left[\left( m({\pmb{\xi}} - {\bf{N}}({\bf{x}})) - m({\pmb{\xi}} -
{\pmb{\xi}}_{0})\right) \right]  \\
\times L_{{\pmb{\xi}}}^{N} \left[\left( \widehat{\phi} \left( \frac{1}{6} |\omega_{+}|^{ - \frac{1}{n}} ({\pmb{\xi}}
- c(\omega_{+})) \right) - \widehat{\phi} \left( \frac{1}{6} |\omega_{-}|^{ - \frac{1}{n}} ({\pmb{\xi}}
- c(\omega_{-})) \right)\right) \right]\, d{\pmb{\xi}} \label{innerint2}
\end{multline}
\begin{multline*}
= \int_{{\pmb{\xi}} \in \mathbb R^{n}} e^{2 \pi i {\bf{z}}
\cdot {\pmb{\xi}}}
\left[\left( m({\pmb{\xi}} - {\bf{N}}({\bf{x}})) - m({\pmb{\xi}} -
{\pmb{\xi}}_{0})\right) \right]  \\
\times \Bigl[ |\omega_{+}|^{- \frac{2N}{n}} \left(\frac{1}{6}\right)^{2N} \left(L_{{\pmb{\xi}}}^{N}
\widehat{\phi}\right) \left( \frac{1}{6} |\omega_{+}|^{ - \frac{1}{n}} ({\pmb{\xi}}
- c(\omega_{+})) \right)  \\ - |\omega_{-}|^{- \frac{2N}{n}} \left(\frac{1}{6}\right)^{2N}
\left(L_{{\pmb{\xi}}}^{N}
\widehat{\phi}\right) \left( \frac{1}{6} |\omega_{-}|^{ - \frac{1}{n}} ({\pmb{\xi}}
- c(\omega_{-})) \right) \Bigr]\, d{\pmb{\xi}}.
\end{multline*}
Therefore the expression in (\ref{innerint2}) is bounded above by a
constant multiple of 
\begin{displaymath}
|\omega_{+}|^{1 - \frac{2N}{n}} + |\omega_{-}|^{1 - \frac{2N}{n}}
 \lesssim |\omega_{-}|^{1 - \frac{2N}{n}}.
\end{displaymath}
All the other terms originating from the integration by parts yield the
same bound. Also note that for ${\bf{z}} \in I_{-} + {\bf{r}}$, we have 
$|{\bf {z}}|^{-2N} \sim \left (|{\mathbf {r}}||I_{-}| \right)^{-2N}$
since $|\mathbf r| \geq 2$. Therefore choosing $N$ large enough, we obtain
\begin{align*}
|\mathcal B_{3}| & \lesssim \sum_{{\bf{r}} \in \mathbb Z^{n} \setminus
 \left\{0 \right\}} \int_{{\bf{z}} \in I_{-} + {\bf{r}}} |{\bf{z}}|^{
 -2N} | \omega_{-}|^{ 1 - \frac{2N}{n}} \left|G_{2}({\bf{x}} - {\bf{z}})\right|\,
 d{\bf{z}}\\
&\lesssim \sum_{{\bf{r}} \in \mathbb Z^{n} \setminus
 \left\{0 \right\}}
| \omega_{-}|^{ 1 - \frac{2N}{n}} \left( |{\bf{r}}|
 |I_{-}|^{\frac{1}{n}} \right)^{-2N} \int_{{\bf{z}} \in I_{-} + \bf{r}}
\left|G_{2}({\bf{x}} - {\bf{z}}) \right|\,
 d{\bf{z}}\\
& \lesssim \sum_{{\bf{r}} \in \mathbb Z^{n} \setminus
 \left\{0 \right\}} |{\bf{r}}|^{n - 2N} \frac{1}{\left(|{\bf{r}}|
 |I_{-}|^{\frac{1}{n}} \right)^{n}} \int_{{\bf{z}} \in C|{\bf{r}}|
I_{-}} \left|G_{2}({\bf{x}} - {\bf{z}}) \right|\,
 d{\bf{z}}\\
& \lesssim \sup_{I \,:\, J \subset I} \frac{1}{|I|} \int_{I}
 |G_{2}({\bf{z}})|\, d{\bf{z}},
\end{align*}
since $J \subset C |{\bf{r}}| I_{-}$. 

It therefore remains to analyze $\mathcal B_{4}$. Since $ {\bf{z}} \in
I_{-} \setminus I_{+}$, there exists 
$$ {\bf{r}} \in \mathbb Z^{n}
\setminus \left\{0\right\}, 1 \leq |{\bf{r}}| \lesssim \left(
\frac{|\omega_{+}|}{|\omega_{-}|} \right)^{\frac{1}{n}}$$
 such that
${\bf{z}} \in I_{+} + {\bf{r}}$. Let us first analyze the integral in ${\bf{y}}$.  Now,
\begin{align*}
&\biggl| \int_{{\bf{y}} \in \mathbb R^{n}} \left( e^{2 \pi i {\bf{N}} \left( {\bf{x}} \right)
\cdot {\bf{y}}} - e^{2 \pi i {\pmb{\xi}}_{0}\cdot{\bf{y}}} \right) K({\bf{y}})
e^{2 \pi i c \left( \omega_{+} \right) \cdot \left({\bf{z}} - {\bf{y}}
\right)}  \\
&  \qquad \qquad \qquad \times \left( \frac{1}{6} | \omega_{+}|^{ - \frac{1}{n}} \right)^{-n} \phi
\left(\frac{{\bf{z}} - {\bf{y}}}{\frac{1}{6} |\omega_{+}|^{- \frac{1}{n}}}
\right) d{\bf{y}} \biggr|                                                          \\
\leq & \sum_{{\bf{m}} \in \mathbb Z^{n}} \biggl|\int_{{\bf{y}} \in I_{+}
+ {\bf{m}}} \left( e^{2 \pi i {\bf{N}} \left( {\bf{x}} \right)
\cdot {\bf{y}}} - e^{2 \pi i {\pmb{\xi}}_{0}\cdot{\bf{y}}} \right) K({\bf{y}})
e^{2 \pi i c \left( \omega_{+} \right) \cdot \left({\bf{z}} - {\bf{y}}
\right)}\\
& \qquad \qquad \qquad \times \left( \frac{1}{6} | \omega_{+}|^{ - \frac{1}{n}} \right)^{-n}  \phi
\left(\frac{{\bf{z}} - {\bf{y}}}{\frac{1}{6} |\omega_{+}|^{- \frac{1}{n}}}
\right)\left( \frac{1}{6} | \omega_{+}|^{ - \frac{1}{n}} \right)^{n} \,  d{\bf{y}} \biggr|                                                              \\ 
\leq & \int_{{\bf{y}} \in I_{+}} |\omega_{-}|^{\frac{1}{n}} |{\bf{y}}|^{1-n}
|\omega_{+}| |{\bf{r}}|^{ - \nu} d{\bf{y}} \quad + \\
&  \left( \sum_{{\bf{m}} \ne 0 \,:\, \left|{\bf{m}} - {\bf{r}}
\right| \leq \frac{|{\bf{r}}|}{2}}\! +\! \sum_{{\bf{m}} \ne 0 \,:\, \left|{\bf{m}} - {\bf{r}}
\right| > \frac{|{\bf{r}}|}{2}} \right) \!\! \left( |{\bf{m}}|
|I_{+}|^{\frac{1}{n}} \right)^{1-n} | \omega_{-}|^{\frac{1}{n}} \min \left(1,
|{\bf{m}} - {\bf{r}}|^{- \nu}| \right) \\
\leq & |\omega_{-}|^{\frac{1}{n}} | \omega_{+} | |I_{+}|^{\frac{1}{n}}
|{\bf{r}}|^{ - \nu} + \left( |{\bf{r}}| |I_{+}|^{\frac{1}{n}} \right)^{1-n} |
\omega_{-} |^{\frac{1}{n}}  \\
& \qquad \qquad \qquad  +\sum_{{\bf{m}} \,:\, \left|{\bf{m}} - {\bf{r}} \right| >
\frac{|{\bf{r}}|}{2}} |{\bf{m}} - {\bf{r}}|^{ - \nu} |\omega_{-}|^{\frac{1}{n}} |I_{+}|^{\frac{1}{n} 
\left( 1 -n\right)}\\
\lesssim & \left( |{\bf{r}}| |I_{+}|^{\frac{1}{n}} \right)^{1-n} |
\omega_{-} |^{\frac{1}{n}} \\
\lesssim & |{\bf{z}}|^{1-n} |\omega_{-}|^{\frac{1}{n}}.
\end{align*}
Therefore, 
\begin{align*}
|\mathcal B_{4}| &\lesssim | \omega_{-}|^{\frac{1}{n}}
 \int_{\left|I_{+}\right|^{\frac{1}{n}} \lesssim |{\bf{z}}| \lesssim
 \left|I_{-}\right|^{\frac{1}{n}}} |{\bf{z}}|^{1-n}
 |G_{2}({\bf{x}} - {\bf{z}})| d{\bf{z}} \\
& \lesssim g \ast |G_{2}|({\bf{x}}),
\end{align*}
where 
\begin{displaymath}
g({\bf{z}}) = h(|{\bf{z}}|) = | \omega_{-} |^{\frac{1}{n}}
|{\bf{z}}|^{1-n} 1_{\left|I_{+}\right|^{\frac{1}{n}} \lesssim |{\bf{z}}| \lesssim
 \left|I_{-}\right|^{\frac{1}{n}}} ({\bf{z}}).
\end{displaymath}
We observe that $ g \in L^{1}, ||g||_{1} \leq C$ and $g$ is radially
decreasing. Let us approximate $g$ from below by $g_{\gamma}$ defined
as follows, 
\[ g_{\gamma} ({\bf{x}}) = 
\begin{cases}
0 &{\mbox{ if }} |{\bf{x}}| \leq |I_{+}|^{\frac{1}{n}} \\
h \left(|I_{+}|^{\frac{1}{n}} + k \gamma \right) 
&{\mbox{ if }} (k-1)\gamma +
|I_{+}|^{\frac{1}{n}} < |{\bf{x}}| \leq k \gamma +
|I_{+}|^{\frac{1}{n}}, \\
& \hspace{4pt} 1 \le k \le k_0\\
0 &{\mbox{ if }} |{\bf{x}}| > |I_{-}|^{\frac{1}{n}}
\end{cases}\]
where $k_0 \gamma = |I_{-}|^{\frac{1}{n}} - |I_{+}|^{\frac{1}{n}}$.\\
We can write $g_{\gamma}$ as 
\begin{multline*}
g_{\gamma} = - h \left(|I_{+}|^{\frac{1}{n}} + \gamma \right) 1_{B
\left(0;\left|I_{+}\right|^{\frac{1}{n}} \right)} \quad  \\ +\sum_{k=1}^{k_0} \left(
h\left(|I_{+}|^{\frac{1}{n}} + k \gamma \right) - h
\left(|I_{+}|^{\frac{1}{n}} + (k+1) \gamma \right) \right)1_{B
\left(0;\left|I_{+}\right|^{\frac{1}{n}} + k \gamma \right)}.
\end{multline*}
Therefore,
\begin{multline*} 
|g_{\gamma} \ast |G_{2}|| \leq  h(|I_{+}|^{\frac{1}{n}} + \gamma)
|G_{2}| \ast 1_{B
\left(0;\left|I_{+}\right|^{\frac{1}{n}} \right)} \quad  \\ +\sum_{k} \left(
h\left(|I_{+}|^{\frac{1}{n}} + k \gamma \right) - h
\left(|I_{+}|^{\frac{1}{n}} + (k+1) \gamma \right) \right) |G_{2}| \ast 1_{B
\left(0;\left|I_{+}\right|^{\frac{1}{n}} + k \gamma \right)},
\end{multline*}
which in turn is bounded by
\begin{align*}
&   \left( \sup_{I \,:\, J \subset I} \frac{1}{|I|} \int_{I}
|G_{2}({\bf{z}})|d{\bf{z}} \right)\;  \times \;\Bigl[h
 \left(|I_{+}|^{\frac{1}{n}} + \gamma \right) \left|{B
\left(0;\left|I_{+}\right|^{\frac{1}{n}} \right)} \right| \quad      \\  
  &  \qquad  +\sum_{k}  \left(
h\left(|I_{+}|^{\frac{1}{n}} + k\gamma \right) - h
\left(|I_{+}|^{\frac{1}{n}} + (k+1) \gamma \right) \right)\left|{B
\left(0;\left|I_{+}\right|^{\frac{1}{n}} + k \gamma \right)}\right| \Bigr]   \\
\leq   &   \left( \sup_{I \,:\, J \subset I} \frac{1}{|I|} \int_{I}
|G_{2}({\bf{z}})|d{\bf{z}} \right)  ||\tilde{g_{\gamma}}||_{1}.
\end{align*}
Here $\tilde{g_{\gamma}}$ is given by 
\begin{multline*}
\tilde{g_{\gamma}} = h \left(|I_{+}|^{\frac{1}{n}} + \gamma \right) 1_{B
\left(0;\left|I_{+}\right|^{\frac{1}{n}} \right)} \quad  \\ +\sum_{k} \left(
h\left(|I_{+}|^{\frac{1}{n}} + k \gamma \right) - h
\left(|I_{+}|^{\frac{1}{n}} + (k+1) \gamma \right) \right)1_{B
\left(0;\left|I_{+}\right|^{\frac{1}{n}} + k \gamma \right)}.
\end{multline*}
In other words, 
\[ \tilde{g_{\gamma}}({\bf{x}}) = 
\begin{cases} 
2 h \left(|I_{+}|^{\frac{1}{n}} + \gamma \right) &{\mbox{ if }} 0 \leq |{\bf{x}}|
\leq |I_{+}|^{\frac{1}{n}} \\
h \left(|I_{+}|^{\frac{1}{n}} + k\gamma \right) &{\mbox{ if }}
|I_{+}|^{\frac{1}{n}} +  (k-1) \gamma < |{\bf{x}}| <
|I_{+}|^{\frac{1}{n}} +  k \gamma\\
0 &{\mbox{ if }} |{\bf{x}}| > |I_{-}|^{\frac{1}{n}}
\end{cases} \]
Therefore,
\begin{displaymath}
||\tilde{g_{\gamma}}||_{1} \lesssim ||g||_{1} +
  |\omega_{-}|^{\frac{1}{n}} \int_{|{\bf{y}}| \leq
|I_{+}|^{\frac{1}{n}}} |{\bf{y}}|^{1-n} d{\bf{y}} \lesssim 1,
\end{displaymath}
since $|\omega_{-}| < | \omega_{+}|$. Thus, 
\begin{displaymath}
g_{\gamma} \ast |G_{2}|({\bf{x}}) \lesssim \sup_{I\,:\, J \subset I}
\frac{1}{|I|} \int_{I}|G_{2}({\bf{z}})| \,  d{\bf{z}}.
\end{displaymath}
Letting $\gamma \rightarrow 0$ and applying the dominated convergence
theorem now yields
the desired bound for $\mathcal B_{4}$. 


\section{An Application}
\setcounter{equation}0

As an immediate application of the weak $L^{2}$ mapping property of the maximal dyadic sum operator, we obtain 
a new proof of Sj\"olin's theorem \cite{sjolin} on a weak-type (2,2) estimate for the maximal conjugated Calder\'on-Zygmund operator on $\mathbb R^{n}$, $n > 1$. 

\begin{theorem}
Let 
\begin{displaymath}
K({\bf{x}}) = \Omega \left(\frac{{\bf{x}}}{|{\bf{x}}|}\right)|{\bf{x}}|^{-n}
\end{displaymath}
be a Calder\'on-Zygmund kernel in $\mathbb R^{n}$ with $ \Omega \in C^{\infty}(S^{n-1})$. Let 
\begin{displaymath}
Bf = f \ast K,
\end{displaymath}
and 
\begin{displaymath}
\mathcal C f({\bf{x}}) = \sup_{{\pmb{\xi}} \in \mathbb R^{n}} \left|\left(e^{2 \pi i
{\pmb{\xi}} \cdot (\cdot)} B e^{- 2 \pi i
{\pmb{\xi}} \cdot (\cdot)}f \right)\right|({\bf{x}}).
\end{displaymath}
Then, 
\begin{displaymath}
||\mathcal C f||_{L^{2,\infty}} \leq C ||f||_{L^{2}},
\end{displaymath}
with a constant $C$ independent of $f$. \label{mainthm}
\end{theorem}

 The proof of Theorem 2 again follows techniques similar to those used by Lacey and Thiele
\cite{lacey-thiele3} in proving Carleson's theorem on almost everywhere convergence of Fourier series. 
Following \cite{lacey-thiele3}, we introduce the operators
\begin{align*}
A_{{\pmb{\eta}}}f & := \sum_{p \in \dd} \langle f, \phi_{p} \rangle
\phi_{p} 1_{\omega_{p(2^{n})}}({\pmb{\eta}})\\
Af & := \lim_{N \rightarrow \infty} \frac{1}{K_{N}} \int_{K_{N} \times
\left[0,1\right]}
M_{- {\pmb{\eta}}} T_{- {\bf{y}}} D_{2^{- \kappa}}^{2} A_{2^{- \kappa} {\pmb{\eta}}}
D_{2^{\kappa}}^{2} T_{{\bf{y}}} M_{{\pmb{\eta}}} f \, d{\bf{y}} \,
d{\pmb{\eta}} \, d\kappa, 
\end{align*}
where $K_{N}$ is any increasing sequence of rectangles filling out
$\mathbb R^{n} \times \mathbb R^{n}$. For any Schwartz function $f$ and
any ${\bf{x}} \in \mathbb R^{n}$, the limit representing
$Af({\bf{x}})$ exists by the argument given by Lacey and Thiele. \\

Note that by rotation invariance, it is enough to prove 
Theorem 2 when the multiplier is supported on a nonempty open cone in
$\mathbb R^{n}$. \\
\begin{lemma}  
There exists a nonempty open cone $ \tilde K_{0}$ with vertex at the origin,
\begin{displaymath}
\tilde K_{0} \subset \left\{{\pmb{\xi}} = (\xi_{1}, \xi_{2}, \cdots, \xi_{n}) ;
\xi_{i} \leq 0 {\mbox{ for all }} i \right\},
\end{displaymath}
such that for all ${\pmb{\xi}} \in \tilde K_{0}$, 
\begin{displaymath}
(Af)\,\widehat{}\,({\pmb{\xi}}) = c {\widehat{f}}({\pmb{\xi}}),
\end{displaymath} 
where $c$ is a constant independent of $f$. 
\end{lemma}
  
\begin{proof}
{\allowdisplaybreaks
\begin{align*}
&(Af)\,\widehat{}\,({\pmb{\xi}}) \\
&= \!\! \lim_{N \rightarrow \infty} \frac{1}{|K_{N}|} 
\intl_{\left( \left({\bf{y}}, {\pmb{\eta}} \right), \kappa \right) \in K_{N} \times
\left[0,1\right]} 
\left(M_{- {\pmb{\eta}}} T_{- {\bf{y}}} D_{2^{- \kappa}}^{2} A_{2^{-
\kappa}  {\pmb{\eta}}}
D_{2^{\kappa}}^{2} T_{{\pmb{y}}} M_{{\pmb{\eta}}} f\right)\widehat{}\,({\pmb{\xi}}) \, d{\bf{y}} \, 
d{\pmb{\eta}} \, d\kappa\\
&= \!\! \lim_{N \rightarrow \infty} \frac{1}{|K_{N}|} \intl_{K_{N} \times
\left[0,1\right]} \sum_{p \in \dd} \langle f, M_{- {\pmb{\eta}}}
T_{-{\bf{y}}} D_{2^{-\kappa}}^{2} \phi_{p} \rangle \! \! \left( M_{-{\pmb{\eta}}} T_{-{\bf{y}}}
D_{2^{-\kappa}}^{2} \phi_{p} \right)\widehat{} \, ({\pmb{\xi}}) \, d{\bf{y}} \, d{\pmb{\eta}} \,
d\kappa\\
& = \!\! \lim_{N \rightarrow \infty}   \frac{1}{|K_{N}|} \intl_{K_{N} \times
\left[0,1\right]} \sum_{p \in \dd}  \int_{\mathbb R^{n}}
\widehat{f}({\pmb{\xi'}}) 2^{\frac{-\kappa n}{2}}
\overline{{\widehat{\phi}_{p}}\left(2^{-\kappa}({{\pmb{\xi'}}}+ {\pmb{\eta}} \right)} e^{- 2 \pi i
{\bf{y}} \cdot \left( {\pmb{\xi'}}+ {\pmb{\eta}} \right)} d {\pmb{\xi'}}      \\ 
& \qquad \qquad \times 2^{\frac{-\kappa n}{2}} \widehat{\phi_{p}}\left(2^{- \kappa} ({\pmb{\xi}} +
{\pmb{\eta}}) \right) e^{2 \pi i {\bf{y}} \cdot \left( {\pmb{\xi}} + {\pmb{\eta}} \right)}
1_{\omega_{p(2^{n})}}(2^{- \kappa }{\pmb{\eta}}) \, d{\bf{y}} \, d{\pmb{\eta}} \, d\kappa \\
&= \!\! \lim_{N \rightarrow \infty}  \frac{1}{|K_{N}|} \intl_{K_{N} \times
\left[0,1\right]} \sum_{p \in \dd}   \int_{\mathbb R^{n}}
\widehat{f}({\pmb{\xi'}}) 2^{- \kappa n +  \alpha n} \overline{{\widehat{\phi}}\left(2^{ \alpha - \kappa}(\xi'+
{\pmb{\eta}}) - \left({\bf{l}}+\frac{1}{4}\right) \right)}  \\ 
& \qquad \qquad \times e^{ - 2 \pi i
\left({\bf{m}}+\frac{1}{2}\right)\cdot
\left(2^{ \alpha - \kappa}\left({\pmb{\xi'}}+{\pmb{\eta}}\right)-\left({\bf{l}}+\frac{1}{4}\right)\right)}
\widehat{\phi}\left(2^{ \alpha -\kappa}({\pmb{\xi}}+
{\pmb{\eta}}) - \left({\bf{l}}+\frac{1}{4}\right) \right)  \\ 
& \qquad \qquad \times e^{ 2 \pi i
\left({\bf{m}}+\frac{1}{2}\right)\cdot
\left(2^{ \alpha -\kappa}\left({\pmb{\xi}} + {\pmb{\eta}}\right)-\left({\bf{l}}+\frac{1}{4}\right)\right)}
 1_{\omega_{p(2^{n})}}(2^{- \kappa }\eta) \, d{\bf{y}} \, d{\pmb{\eta}}  \, d\kappa \, d{\pmb{\xi'}},
\end{align*}}
where we have expressed $I_{p}$ and $\omega_{p}$ as
\begin{displaymath}
I_{p} = \prod_{i=1}^{n} \left[m_{i}2^{\alpha}, (m_{i}+1)2^{\alpha}\right), \quad 
\omega_{p} = \prod_{i=1}^{n} \left[l_{i}2^{-\alpha }, (l_{i}+1)2^{- \alpha}\right)
\end{displaymath}
with 
\begin{displaymath}
{\bf{m}} = (m_{1}, m_{2}, \cdots , m_{n}), \;  {\bf{l}} = (l_{1}, l_{2},
\cdots  l_{n}) \in \mathbb Z^{n}.
\end{displaymath}
Interpreting the sum 
\begin{displaymath}
\sum_{{\bf{m}} \in \mathbb Z^{n}} e^{2 \pi i \left({\bf{m}} + \frac{1}{2}\right)\cdot
2^{ \alpha - \kappa}\left( {\pmb{\xi}} - {\pmb{\xi'}} \right)}
\end{displaymath}
in the sense of distributions, we find that 
\begin{multline}
(Af)\,\widehat{}\,({\pmb{\xi}}) = c \lim_{N \rightarrow \infty}\frac{1}{|K_{N}|}
\int_{K_{N} \times \left[0,1\right]} \sum_{ {\bf{l}} \in \mathbb Z^{n},
\alpha  \in \mathbb Z} \widehat{f}({\pmb{\xi}}) \\
\times \left|{\widehat{\phi}}\left(2^{ \alpha -\kappa} \left({\pmb{\xi}} +
{\pmb{\eta}} \right) -  \left( {\bf{l}} + \frac{1}{4} \right) \right)
\right|^{2} 1_{\omega_{p(2^{n})}}(2^{-\kappa} {\pmb{\eta}})  \,
d{\bf{y}} \, d{\pmb{\eta}} \,
d\kappa. \label{af}
\end{multline}
Now, ${\widehat{\phi}}\left(2^{ \alpha -\kappa} \left({\pmb{\xi}} +
{\pmb{\eta}} \right) -  \left( {\bf{l}} + \frac{1}{4} \right) \right)$ is
supported on 
\begin{multline*}
\biggl\{({\pmb{\xi}}, {\pmb{\eta}}) : -\frac{1}{10} \leq 2^{ \alpha -\kappa} \left(
\xi_{i} + \eta_{i} \right) - \left( l_{i} + \frac{1}{4} \right) \leq
\frac{1}{10} {\mbox{ for all }} i, \, 1 \leq i \leq n \biggr\}=\\
\biggl\{({\pmb{\xi}}, {\pmb{\eta}}) : \left(l_{i} + \frac{3}{20}
\right) 2^{\kappa- \alpha } \leq \xi_{i} + \eta_{i} \leq \left( l_{i} + \frac{7}{20}
\right) 2^{\kappa - \alpha } {\mbox{ for all }} i, \, 1 \leq i \leq n \biggr\}.
\end{multline*}
Also, 
\begin{align*}
1_{\omega_{p(2^{n})}}(2^{- \kappa} {\pmb{\eta}}) \ne 0 & \Leftrightarrow
\left(l_{i}+ \frac{1}{2} \right)2^{-p} \leq 2^{- \kappa} \eta_{i} \leq
\left(l_{i} + 1 \right) 2^{- \alpha } \\
& \Leftrightarrow \left(l_{i}+ \frac{1}{2} \right)2^{\kappa - \alpha } \leq \eta_{i} \leq
\left(l_{i} + 1 \right) 2^{\kappa - \alpha }.
\end{align*}
Therefore, the integrand in the right hand side of (\ref{af}) is supported in 
\begin{displaymath}
\bigcup_{ \alpha \in \mathbb Z}\left\{ - \frac{17}{20} 2^{\kappa - \alpha} \leq \xi_{i} \leq - \frac{3}{20}
2^{\kappa - \alpha} {\mbox{ for all }} i,\, 1 \leq i \leq n  \right\}.
\end{displaymath}
Moreover, if \begin{equation}
- \frac{13}{20} 2^{\kappa - \alpha } \leq \xi_{i} \leq - \frac{7}{20}
2^{\kappa - \alpha} {\mbox{ for all }} i, \quad 1 \leq i \leq n, \label{cone}
\end{equation}
then 
\begin{multline*}
\left\{{\pmb{\eta}} \,:\, \left(l_{i} + \frac{3}{20}
\right) 2^{\kappa - \alpha} \leq \xi_{i} + \eta_{i} \leq \left( l_{i} + \frac{7}{20}
\right) 2^{\kappa - \alpha} {\mbox{ for all }} i, \, 1 \leq i \leq n \right\}\\
\subset  \left\{ {\pmb{\eta}} \, : \,\left(l_{i}+ \frac{1}{2}
\right)2^{\kappa - \alpha} \leq \eta_{i} \leq
\left(l_{i} + 1 \right) 2^{\kappa - \alpha} \right\}. 
\end{multline*}
Note that there exists a nonempty open cone $ \tilde K_{0}$ with vertex at the origin such
that for all ${\pmb{\xi}} \in \tilde{K_{0}}$,
there exist $\alpha \in \mathbb Z$ and $ \kappa \in [0,1]$ satisfying
(\ref{cone}). In the sequel we work with such a cone. 

For ${\pmb{\xi}} \in \tilde K_{0}$ and choosing $K_{N} = [-N,N]^{n} \times [-N, N]^{n}$, we have
{\allowdisplaybreaks \begin{align*}
&(Af)\,\widehat{}\,({\pmb{\xi}}) \\
& = c \widehat{f}({\pmb{\xi}})\; 
 \lim_{N \rightarrow \infty}\frac{1}{|K_{N}|}\\
 & \times
\int_{K_{N} \times \left[0,1\right]} \sum_{ {\bf{l}} \in \mathbb Z^{n},
\alpha \in \mathbb Z} \left|{\widehat{\phi}}\left(2^{\alpha -\kappa} \left({\pmb{\xi}} +
{\pmb{\eta}} \right) -  \left( {\bf{l}} + \frac{1}{4} \right) \right)
\right|^{2}\, d{\bf{y}} \, d{\pmb{\eta}} \,d\kappa \\
& = c \widehat{f}({\pmb{\xi}}) \; 
 \lim_{N \rightarrow \infty}\frac{1}{N^{n}} \\
 & \times
\int_{\left({\pmb{\eta}}, \kappa \right) \in \left[-N, N\right]^{n}
\times \left[0,1\right]} \!\! \sum_{\begin{subarray}{c}{\bf{l}} \in \mathbb Z^{n},
p \in \mathbb Z\\ \left|\xi_{i}\right| \sim 2^{\kappa - \alpha} \\ \eta_{i} \sim l_{i}
2^{\kappa - \alpha} \, \forall i \end{subarray}} \!\!
\left|{\widehat{\phi}}\left(2^{\alpha -\kappa} \left({\pmb{\xi}} +
{\pmb{\eta}} \right) -  \left( {\bf{l}} + \frac{1}{4} \right) \right)
\right|^{2} d{\pmb{\eta}} \,d\kappa \\
&= c \widehat{f}({\pmb{\xi}}) \lim_{N \rightarrow \infty} N^{-n}
\int_{0}^{1} \sum_{\begin{subarray}{c}{\bf{l}} \in \mathbb Z^{n}, \alpha \in
\mathbb Z\\\left|\xi_{i}\right| \sim 2^{\kappa -\alpha}\\ -N \lesssim l_{i} 2^{ \kappa
- \alpha} \lesssim N \, \forall i \end{subarray}} \left(2^{\kappa - \alpha} \right)^{n} \, d\kappa \\
&= c \widehat{f}({\pmb{\xi}}) \lim_{N \rightarrow \infty} N^{-n}
\int_{0}^{1}\sum_{\begin{subarray}{c}\alpha \in
\mathbb Z\\\left|\xi_{i}\right| \sim 2^{\kappa - \alpha}\, \forall i \end{subarray}}
\left(N2^{-\kappa + \alpha}\right)^{n} \left(2^{\kappa- \alpha}\right)^{n} d\kappa \\
&= c \widehat{f}({\pmb{\xi}})\int_{0}^{1} \sum_{\begin{subarray}{c}\alpha \in
\mathbb Z\\\left|\xi_{i}\right| \sim 2^{\kappa - \alpha}\, \forall i \end{subarray}} 1 \, 
d\kappa \\
&= c \widehat{f}({\pmb{\xi}}).
\end{align*}}
\end{proof}

Now, let $m$ be the multiplier associated with the Calder\'on-Zygmund
kernel $K$. Then $m \in C^{\infty}(\mathbb R^{n} \setminus
\{\mathbf 0\})$ and is homogeneous of degree 0. Suppose further, without loss of
generality, that $m$ is supported on the cone $ \tilde K_{0}$ described
earlier. We may reduce the problem to this case via a partition of unity and by
invoking rotation invariance. \\
\noindent Recalling that 
\begin{displaymath}
Bf = f \ast K,
\end{displaymath}
and
\begin{displaymath}
\mathcal C f({\bf{x}}) = \sup_{{\pmb{\eta}}} |M_{{\pmb{\eta}}} B M_{-{\pmb{\eta}}}f|({\bf{x}}),
\end{displaymath}
we define, for ${\pmb{\zeta}} \in \mathbb R^{n}$, 
\begin{displaymath}
B_{{\pmb{\zeta}}}f := \sum_{p \in \dd} \langle f, \phi_{p}
\rangle \psi_{p}^{{\pmb{\zeta}}} 1_{\omega_{p(2^{n})}}({\pmb{\zeta}}),
\end{displaymath}
where 
\begin{displaymath}
\left( \psi_{p}^{{\pmb{\zeta}}} \right)\, \widehat{} \, ({\pmb{\xi}}) = m( {\pmb{\xi}} -
{\pmb{\zeta}}) \widehat{\phi_{p}}({\pmb{\xi}}).
\end{displaymath}

Using Lemma 6, it is not hard to see that
\begin{align*}   
Bf :&= \lim_{N \rightarrow \infty} \frac{1}{|K_{N}|} \int_{K_{N} \times
\left[0,1\right]}
M_{- {\pmb{\eta}}} T_{- {\bf{y}}} D_{2^{- \kappa}}^{2} B_{2^{- \kappa} {\pmb{\eta}}}
D_{2^{\kappa}}^{2} T_{{\bf{y}}} M_{{\pmb{\eta}}} f \, d{\bf{y}} \,
d{\pmb{\eta}} \, d\kappa\\ 
& = \lim_{N \rightarrow \infty} \frac{1}{|K_{N}|} \int_{K_{N} \times
\left[0,1\right]} \sum_{p \in \dd} \Bigl[\langle f, M_{- {\pmb{\eta}}}
T_{-{\bf{y}}} D_{2^{-\kappa}}^{2} \phi_{p} \rangle \;  \\
& \qquad \qquad \qquad \times \left(M_{-{\pmb{\eta}}} T_{-{\bf{y}}}
D_{2^{-\kappa}}^{2} \psi_{p}^{2^{-\kappa} {\pmb{\eta}}}\right) \,
\left(1_{\omega_{p(2^{n})}}(2^{- \kappa} {\pmb{\eta}})\right) \Bigr] \,
d{\bf {y}} \, d{\pmb{\eta}} \,d\kappa. 
\end{align*}
In fact, 
\begin{align*}
(Bf)\,\widehat{}\, ({\pmb{\xi}}) &= \lim_{N \rightarrow \infty} \frac{1}{|K_{N}|} \int_{K_{N} \times
\left[0,1\right]} \sum_{p \in \dd}  \Bigl[\langle f, M_{- {\pmb{\eta}}}
T_{-{\bf{y}}} D_{2^{-\kappa}}^{2} \phi_{p} \rangle \;  \\
& \qquad  \times T_{-{\pmb{\eta}}} M_{{\bf{y}}}
D_{2^{\kappa}}^{2} \left(\psi_{p}^{2^{-\kappa} {\pmb{\eta}}} \right)
\widehat{} \left( {\pmb{\xi}} \right) \,
\left(1_{\omega_{p(2^{n})}}(2^{- \kappa} {\pmb{\eta}})\right) \Bigr] \,
d{\bf {y}} \, d{\pmb{\eta}} \,d\kappa \\
&  = \lim_{N \rightarrow \infty} \frac{1}{|K_{N}|} \int_{K_{N} \times
\left[0,1\right]} \sum_{p \in \dd}  \Bigl[\langle f, M_{- {\pmb{\eta}}}
T_{-{\bf{y}}} D_{2^{-\kappa}}^{2} \phi_{p} \rangle \; \\
& \qquad  \times T_{-{\pmb{\eta}}} M_{{\bf{y}}}
D_{2^{\kappa}}^{2} \left[m \left({\pmb{\xi}} - {\pmb{\eta}} 2^{- \kappa}\right)
\widehat{\phi_{p}}({\pmb{\xi}}) \right] \,
\left(1_{\omega_{p(2^{n})}}(2^{- \kappa} {\pmb{\eta}})\right) \Bigr] \,
d{\bf {y}} \, d{\pmb{\eta}} \,d\kappa. 
\end{align*}
Since
\begin{displaymath}
T_{-{\pmb{\eta}}} M_{{\bf{y}}}
D_{2^{\kappa}}^{2}
\left[m \left({\pmb{\xi}} - {\pmb{\eta}} 2^{- \kappa}\right)
\widehat{\phi_{p}}({\pmb{\xi}}) \right] = m({\pmb{\xi}}) T_{-{\pmb{\eta}}} M_{{\bf{y}}}
D_{2^{\kappa}}^{2}\widehat{\phi_{p}}({\pmb{\xi}}),
\end{displaymath}
we get that
\begin{displaymath}
(Bf)\,\widehat{}\, ({\pmb{\xi}}) = m({\pmb{\xi}}) (Af)\,\widehat{}\,
({\pmb{\xi}}) = c\; m({\pmb{\xi}}) \widehat{f}({\pmb{\xi}}),
\end{displaymath}
where the last equality follows from the claim and the fact that ${\mbox{supp}}\, m
\subset \tilde K_{0}$.

Now, for $f \in C_{0}^{\infty}(\mathbb R^{n})$ and ${\pmb{\zeta}} \in \rn$
\begin{align*}
& M_{{\pmb{\zeta}}}BM_{-{\pmb{\zeta}}}f({\bf{x}})\\
&= \lim_{l \rightarrow \infty} \frac{1}{\left|
K_{l}\right|} \int_{K_{l} \times \left[0,1 \right]} M_{{\pmb{\zeta}}}M_{-{\pmb{\eta}}}
T_{-{\bf{y}}} D_{2^{-\kappa}}^{2} B_{2^{-\kappa}{\pmb{\eta}}} D_{2^{\kappa}}^{2}
T_{\bf{y}} M_{\pmb{\eta}} M_{-{\pmb{\zeta}}}f({\bf{x}}) \, d{\bf{y}}\, d{\pmb{\eta}}
\,d\kappa \\
& = \lim_{l \rightarrow \infty} \frac{1}{\left|
K_{l}\right|} \int_{K_{l} \times \left[0,1 \right]} M_{{\pmb{\zeta}}-\pmb{\eta}}
T_{-\bf{y}} D_{2^{-\kappa}}^{2} B_{2^{-\kappa}{\pmb{\eta}}} D_{2^{\kappa}}^{2}
T_{\bf{y}} M_{\pmb{\eta}-{\pmb{\zeta}}}f({\bf{x}}) \, d{\bf{y}} \,d{\pmb{\eta}}
\,d\kappa \\
&=\lim_{l \rightarrow \infty} \frac{1}{\left|
K_{l}\right|} \int_{K_{l}^{{\pmb{\zeta}}} \times \left[0,1 \right]} M_{-\pmb{\eta'}}
T_{-\bf{y}} D_{2^{-\kappa}}^{2} B_{2^{-\kappa}\left({\pmb{\eta'}} + {\pmb{\zeta}}\right)} D_{2^{\kappa}}^{2}
T_{\bf{y}} M_{\pmb{\eta'}}f({\bf{x}}) \, d{\bf{y}} \, d{\pmb{\eta'}}
\, d\kappa \\
&=\lim_{l \rightarrow \infty} \frac{1}{\left|
K_{l}\right|} \int_{K_{l} \times \left[0,1 \right]} M_{-\pmb{\eta'}}
T_{-\bf{y}} D_{2^{-\kappa}}^{2} B_{2^{-\kappa}\left({\pmb{\eta'}} + {\pmb{\zeta}}\right)} D_{2^{\kappa}}^{2}
T_{\bf{y}} M_{\pmb{\eta'}}f({\bf{x}}) \, d{\bf{y}} \, d{\pmb{\eta'}}
\,d\kappa,
\end{align*}
where $K_{l}^{{\pmb{\zeta}}} := {{\pmb{\zeta}}} + K_{l}$. The last equality follows from
the fact that for $f \in C_{0}^{\infty}(\mathbb R^{n})$ and any fixed
${\pmb{\zeta}}$, the integrand
gets arbitrarily small on the domain $K_{l}^{{\pmb{\zeta}}} \triangle \, K^{l}$ as
$l \rightarrow \infty$. The details are left to the interested reader.

Therefore,
\begin{align*}
 &\sup_{{\pmb{\zeta}}}\left| M_{{\pmb{\zeta}}}BM_{-{\pmb{\zeta}}}f(\bf{x})\right|\\
  &\leq \lim_{l \rightarrow \infty} \frac{1}{\left|
K_{l}\right|} \int_{K_{l} \times \left[0,1 \right]} \left|M_{-\pmb{\eta'}}
T_{-\bf{y}} D_{2^{-\kappa}}^{2} \left[\sup_{{\pmb{\zeta}}}B_{{\pmb{\zeta}}}\right] D_{2^{\kappa}}^{2}
T_{\bf{y}} M_{\pmb{\eta'}}f(\bf{x}) \right|\, d{\bf{y}} \, d{\pmb{\eta'}}
\, d\kappa.
\end{align*}
We recall the following fact about the weak $L^2$ norm :
 there exist universal constants $C_1, C_2$ such that
\[ C_1 \sup_{E}
\frac{\left| \langle f, {{1}}_{E}\rangle \right|}{|E|^{\frac{1}{2}}}
\leq \left|\left| f\right|\right|_{L^{2, \infty}} \leq C_2 \sup_{E}
\frac{\left| \langle f, {{1}}_{E}\rangle \right|}{|E|^{\frac{1}{2}}},
\]
where the supremum is over all measurable sets $E$ with finite Lebesgue
measure. This implies
{\allowdisplaybreaks
\begin{align*}
&\left|\left| \sup_{{\pmb{\zeta}}}\left|
M_{{\pmb{\zeta}}}BM_{-{\pmb{\zeta}}}f \right|\right|\right|_{L^{2,\infty}}\\
 &\lesssim
 \sup_{E} \frac{1}{|E|^{\frac{1}{2}}}
\int_{E} \biggl[\lim_{l \rightarrow \infty} \frac{1}{\left|
K_{l}\right|} \\
& \qq \times \int_{K_{l} \times \left[0,1 \right]} \biggl|M_{-\pmb{\eta'}}
T_{-\bf{y}} D_{2^{-\kappa}}^{2} \biggl(\sup_{{\pmb{\zeta}}}B_{{\pmb{\zeta}}}\biggr) D_{2^{\kappa}}^{2}
T_{\bf{y}} M_{\pmb{\eta'}}f(\bf{x}) \biggr|\, d{\bf{y}} d{\pmb{\eta'}}
d\kappa \biggr]\, d{\bf{x}}\\
& \lesssim  \lim_{l \rightarrow \infty} \frac{1}{\left|
K_{l}\right|}
 \int_{K_{l} \times \left[0,1 \right]} \biggl[\sup_{E} \frac{1}{|E|^{\frac{1}{2}}}\\
 & \qq \times \int_{E} \biggl|M_{-\pmb{\eta'}}
T_{-\bf{y}} D_{2^{-\kappa}}^{2} \biggl(\sup_{{\pmb{\zeta}}}B_{{\pmb{\zeta}}}\biggr) D_{2^{\kappa}}^{2}
T_{\bf{y}} M_{\pmb{\eta'}}f(\bf{x}) \biggr|\, d{\bf{x}} \biggr]
d{\bf{y}}\, d{\pmb{\eta'}} \, d\kappa \\
& \lesssim  \lim_{l \rightarrow \infty} \frac{1}{\left| K_{l}\right|}
\int_{K_l \times \left[0,1 \right]} \left|\left| |M_{-\pmb{\eta'}}
T_{-\bf{y}} D_{2^{-\kappa}}^{2} \left(\sup_{{\pmb{\zeta}}}B_{{\pmb{\zeta}}}\right) D_{2^{\kappa}}^{2}
T_{\bf{y}} M_{\pmb{\eta'}}f |\right|\right|_{L^{2,
\infty}} \!\!\! d{\bf{y}}\, d{\pmb{\eta'}}\, d{\kappa}\\
& \lesssim  \lim_{l \rightarrow \infty} \frac{1}{\left| K_{l}\right|}
\int_{K_l \times \left[0,1 \right]} \left|\left|\left(\sup_{{\pmb{\zeta}}}B_{{\pmb{\zeta}}}\right) D_{2^{\kappa}}^{2}
T_{\bf{y}} M_{\pmb{\eta'}}f \right|\right|_{L^{2,
\infty}} \, d{\bf{y}}\, d{\pmb{\eta'}}\, d{\kappa},
\end{align*}}
since the weak $L^2$ norm is invariant under the translation, dilation and
modulation operators defined in Section 2.  The same invariance properties
also hold true for  the $L^2$ norm.
In order to prove the  weak $L^{2}$ bound for the Carleson
operator, it therefore suffices to show that
\begin{equation}
 || \sup_{{\pmb{\zeta}}} B_{{\pmb{\zeta}}}f ||_{L^{2, \infty}} \leq C ||f||_{L^{2}},
    \label{mainest}
\end{equation} 
which is the conclusion of Theorem 1.

\section{Acknowledgements}
\setcounter{equation}0

Versions of the main result in this paper were obtained independently by the two authors as
part of their dissertations at the University of California, Berkeley and the University of Missouri,
Columbia respectively.  The authors would like to thank their advisors
Michael Christ and Loukas Grafakos for encouragement and many insightful discussions during
preparation of the paper.

\end{document}